\documentclass[11pt]{amsart}
\usepackage{array, amsmath, amscd, amssymb, amsthm}
\usepackage{hyperref}

\newcommand{\inclu}{\hookrightarrow}
\newcommand{\ra}{\longrightarrow}
\newcommand{\ptp}{\PP^1\times \PP^1}
\newcommand{\SL}{\omega_{2,1}}
\newcommand{\pin}{\pi^{[n]}}
\newcommand{\cXn}{\cX^{[n]}}

\newcommand{\Ln}{\cL^{[n]}}
\newcommand{\cLn}{\cL^{[n]}}

\newcommand{\qqq}{\QQ[[q]]^{\times}}
\newcommand{\qqx}{\QQ[[x]]^{\times}}
\newcommand{\Gott}{G\"{o}ttsche}
\newcommand{\Pan}{Pandharipande}
\newcommand{\Pand}{Pandharipande}
\newcommand{\bl}{\mathrm{bl}}
\newcommand{\pr}{\mathrm{pr}}
\newcommand{\fourtop}{L^2, LK, c_1(S)^2, c_2(S)}
\newcommand{\fourtopand}{$L^2$, $LK$, $c_1(S)^2$ and $c_2(S)$}
\renewcommand{\deg}{\mathrm{deg}}
\newcommand{\ds}{\displaystyle}

\theoremstyle{plain}

\newtheorem{prop}{Proposition}[section]
\newtheorem{theo}[prop]{Theorem}
\newtheorem{lemm}[prop]{Lemma}
\newtheorem{coro}[prop]{Corollary}

\newtheorem*{univpoly}{Theorem 1.1}
\newtheorem*{GYZ}{Theorem 1.2}
\newtheorem*{power}{Theorem 1.3}

\theoremstyle{definition}
\newtheorem{defn}{Definition}[section]
\newtheorem{depr}[prop]{Definition-Proposition}

\theoremstyle{remark}
\newtheorem*{rem}{Remark}
\newtheorem*{notation}{Notation}

\newcommand{\CC}{\mathbb{C}}

\newcommand{\FF}{\mathbb{F}}
\newcommand{\LL}{\mathbb{L}}

\newcommand{\NN}{\mathbb{N}}
\newcommand{\PP}{\mathbb{P}}
\newcommand{\QQ}{\mathbb{Q}}
\newcommand{\RR}{\mathbb{R}}

\newcommand{\ZZ}{\mathbb{Z}}

\def\cF{{\mathcal F}}

\def\cI{{\mathcal I}}

\def\cL{{\mathcal L}}
\def\cM{{\mathcal M}}

\def\cO{{\mathcal O}}

\def\cR{{\mathcal R}}

\def\cV{{\mathcal V}}
\def\cW{{\mathcal W}}
\def\cX{{\mathcal X}}

\def\cZ{{\mathcal Z}}

\def\fB{\mathfrak{B}}

\def\fI{\mathfrak{I}}

\let\oldTitle\title
\renewcommand{\title}[1]{\newcommand{\myTitle}{#1}\oldTitle{#1}} 
\title{A Proof of the G\"{o}ttsche-Yau-Zaslow Formula}

\input xy
\xyoption {all}
\numberwithin{equation}{section}
\begin{document}

\author{Yu-jong Tzeng}
\address{Department of Mathematics, Harvard University,
Cambridge, MA 02138, USA} \email{ytzeng@math.harvard.edu}

\begin{abstract}  
Let $S$ be a complex smooth projective surface and $L$ be a line bundle on $S$. 
\Gott\ conjectured that for every integer $r$, the number of $r$-nodal curves in $|L|$ is a universal polynomial of four topological numbers when $L$ is sufficiently ample. 
We prove \Gott's conjecture using the algebraic cobordism group of line bundles on surfaces and degeneration of Hilbert schemes of points.
In addition, we prove the the G\"{o}ttsche-Yau-Zaslow Formula which expresses the generating function of the numbers of nodal curves in terms of quasi-modular forms and two unknown series. 
\end{abstract}
\maketitle

\section{Introduction}\label{intro}

\subsection{Main results}
Consider a line bundle $L$ on a complex projective smooth  surface $S$. 
This paper attempts to answer the following question:  how many reduced curves have exactly $r$ simple nodes and no higher singularities in a generic $r$-dimensional linear system of $|L|$?
Equivalently, how many $r$-nodal curves in $|L|$ pass through dim$|L|-r$ points in general position? 

\Gott\ \cite{Gott} conjectured that for every $r$, the numbers of $r$-nodal curves are given by universal polynomials of four topological numbers: \fourtopand\, provided that the line bundle $L$ is $(5r-1)$-very ample (i.e. $H^0(S,L)\to H^0(L|_{\xi})$ is surjective for every $\xi\in S^{[5r]}$). 
These polynomials are universal in the sense that they only depend on $r$ and are independent of the surface and line bundle. Our first result is an algebro-geometric proof of \Gott's conjecture.
   
\begin{theo}[\Gott's conjecture]\label{univpoly}For every integer $r\geq 0$, there exists a universal polynomial $T_r(x,y,z,t)$ of degree $r$ with the following property: given a smooth projective surface $S$ and a $(5r-1)$-very ample (5-very ample if $r= 1$) line bundle $L$ on $S$, a general $r$-dimensional sublinear system of $|L|$ contains exactly $T_r(L^2, LK, c_1(S)^2, c_2(S))$ $r$-nodal curves.  \end{theo}

Since the numbers of nodal curves for all line bundles on $\PP^2$ and for primitive classes on K3 surfaces have been determined (\cite{CH}, \cite{BL}), all coefficients of $T_r$ can be computed by solving linear equations. 
Moreover, one can combine all universal polynomials as coefficients to define a generating function. 
Inspired by the Yau-Zaslow formula, \Gott\ \cite{Gott} conjectured the closed form of this generating function, which we call the \Gott-Yau-Zaslow formula. 

\begin{theo}[The \Gott-Yau-Zaslow formula]\label{DG2}
There exist universal power series $B_1(q)$ and $B_2(q)$  such that
$$\sum_{r\geq 0}T_r(L^2, LK, c_1(S)^2, c_2(S))(DG_2(\tau))^r=
\frac{(DG_2(\tau)/q)^{\chi(L)}B_1(q)^{K_S^2}B_2(q)^{LK_S}}{(\Delta(\tau)D^2G_2(\tau)/q^2)^{\chi(\cO_S)/2}},$$
where $G_2$ is the second Eisenstein series $-\frac{1}{24}+ \sum_{n>0}\left(\sum_{d|n} d\right)q^n$, $D=q\frac{d}{d\,q}$ and
$\Delta(q)=q\prod_{k>0}(1-q^k)^{24}$. 
\end{theo}

If we write $q=e^{2\pi i\tau}$ for $\tau$ on the complex upper half plane, then $G_2$, $DG_2$ and $D^2G_2$ are quasimodular forms in $\tau$ and $\Delta$ is a modular form in $\tau$. For the precise definition of quasimodular forms, see \cite{KZ}. 

\subsection{Background}\label{history}
On $\PP^2$, the number of nodal curves is classically known as the Severi degree $N_{d,g}$, which is the number of plane curves of degree $d$ and genus $g$ passing through $3d + g - 1$ points in general position. 
This subject was studied by Ran \cite{Ran1} and \cite{Ran2}, using degeneration of $\PP^2$ and an inductive procedure. 
In 1993, Kontsevich and Manin \cite{KM} introduced the techniques of  Gromov-Witten theory to this problem from which they obtained a beautiful recursive formula of rational curves for all degrees. 
For plane curves of a higher genus, Harris and Pandharipande \cite{HP} computed the Severi degrees with at most three nodes using Hilbert schemes, and Choi \cite{Choi} extended the result to at most four nodes using Ran's method.  

The counting of nodal curves of arbitrary genus $g$ in $\PP^2$ was completely solved by Caporaso and Harris \cite{CH}. They defined the generalized Severi degrees $N_{d,g}(\alpha, \beta)$ with tangential conditions and used deformation theory to derive recursive formulas of $N_{d,g}(\alpha, \beta)$. 
Shortly after, Vakil \cite{Vakil} derived similar results for rational ruled surfaces. 

For an arbitrary smooth projective surface, the number of nodal curves with at most three nodes can be computed directly by standard intersection theory. 
In 1994, Vainsencher \cite{Vain} proved the existence of universal polynomials in the case of up to six nodes.
By computing the polynomials explicitly, he showed that the polynomials only depend on $L^2$, $LK$, $c_1(S)^2$ and $c_2(S)$. 
Later, Kleiman and Piene \cite{KP} refined Vainsencher's approach and  generalized the result to up to eight nodes. 
Since their methods rely on a detailed analysis of the singularities of low codimensions, it is difficult to generalize the methods to the case of higher number of nodes. 

On algebraic K3 surfaces and primitive classes, the number of rational curves not only can be determined but also possess a general pattern. 
Yau and Zaslow \cite{YZ} discovered a surprising formula for the generating function in terms of the Dedekind function, which prompts the speculation that general modular forms may be involved.
In particular, the Yau-Zaslow formula implies that the number of rational curves in an effective class $C$ only depends on the self-intersection number $C^2$.

The Yau-Zaslow formula was generalized by \Gott\ \cite{Gott} to arbitrary projective surface. Using Vainsencher and Kleiman-Piene's numbers, \Gott\ conjectured the \Gott-Yau-Zaslow formula which relates the generating function to quasi-modular forms. 
This generating function is defined by universal polynomials, which can be viewed as a virtual counting of nodal curves especially when the line bundle is not ample enough. 
In addition, \Gott\ observed that the generating function is completely determined by the number of nodal curves on $\PP^2$ and K3 surfaces; thus it can be computed by the Severi degrees and quasi-modular forms.  A reformulation of the \Gott-Yau-Zaslow formula gives the generating function on the number of genus $g$ curves (\cite{Gott} Remark 2.6) and this reformulation has been verified by  Bryan and Leung \cite{BL} for  K3 surfaces and primitive ample line bundles. 

The Severi degrees also have interesting properties. On $\PP^2$, fixing the number of nodes and letting the degrees vary, Di Francesco and Itzykson conjectured \cite{FI} that the number of plane $r$-nodal curves of degrees $d$ is a polynomial in $d$, or the node polynomial.  
Recently, Fomin and Mikhalkin \cite{FM} proved the polynomiality with tropical geometry and find many interesting properties of node polynomials. Block \cite{Bl} generalized it to relative node polynomials and proved that there is a formal power series which specializes to all relative node polynomials. These results suggest that enumerating curves with broader conditions may possess a generalized structure, which could be used to provide answers to open problems and interpretations for known results.

We add that a symplectic proof to \Gott's conjecture was given by A.K. Liu \cite{Liu1}, \cite{Liu8}, based on the work of Taubes on the equivalence of Seiberg-Witten theory and Gromov-Witten theory. Recently another proof is also found by Kool-Shende-Thomas \cite{KST}, using the BPS calculus and the computation of tautological integrals on Hilbert schemes by Ellingsrud, \Gott\ and Lehn.

\subsection{Multiplicative Structure}\label{GYZ}
The \Gott-Yau-Zaslow formula (Theorem  \ref{DG2}) expressed the generating function in terms of quasimodular forms and two unknown series $B_1(q)$ and $B_2(q)$.
If the canonical divisor of a surface $S$ is numerically trivial, only the quasimodular forms appear in the generating function.
Consequently, the generating function for K3 surfaces and abelian surfaces is known.   
The unknown series $B_1(q)$ and $B_2(q)$ can be determined using Caporaso and Harris' recursive formulas of the Severi degrees of $\PP^2$. 
G\"{o}ttsche has computed the coefficients of $B_1(q)$ and $B_2(q)$ up to degree 28. 
However, we are still unable to find the closed forms for $B_1(q)$ and $B_2(q)$. 

Instead of indexing on $(DG_2)^r$, we can simply use $x^r$ to define another generating function
$$T(S,L)=\sum_{r\geq0} T_r(L^2, LK, c_1(S)^2, c_2(S))x^r.$$
Because all coefficients are universal, $T(S,L)$ is a universal power series. 
Moreover, $T(S,L)$ is multiplicative: 
 
\begin{theo}[\cite{Gott}, Proposition 2.3]  \label{power} Assuming the numbers of nodal curves are given by universal polynomials, then there exist universal power series $A_1$, $A_2$, $A_3$, $A_4$  in $\qqx$ \footnote{$\qqx$ is the group of units in $\QQ[[x]]$ and the group action is defined by multiplication of power series. } such that the generating function has the form

$$T(S,L)=A_1^{L^2}A_2^{LK_S}A_3^{c_1(S)^2}A_4^{c_2(S)}.$$  \end{theo}

The coefficients of $A_i$ can be determined by Caporaso-Harris \cite{CH} and Vakil's \cite{Vakil} recursive formulas on $\PP^2$ and $\ptp$ but the closed forms are unknown. 

While \Gott\ proved this theorem by considering disjoint union of surfaces, we will give a different proof using algebraic cobordism in Section \ref{chap4}. The new proof plays a central role in our approach because it demonstrates that universality is a result of algebraic cobordism structure.

The difficulty in proving closed formulas of generating functions comes from the fact that they are defined by universal polynomials, not the actual number of nodal curves. The generating functions in Theorem \ref{DG2} and Theorem \ref{power} is well-defined for all line bundles even  if they are trivial or negative.
In this setting, the $r$-th coefficient of $T(S,L)$ equals the number of $r$-nodal curves of $[S,L]$ if $L$ is sufficiently ample relative to $r$. Therefore usually only a finite number of initial coefficients honestly represent the number of nodal curves, and after that the coefficients lack a geometric interpretation. 
This difficulty will be overcome by depicting the universal polynomials as intersection numbers $d_r(S,L)$ and work with $d_r(S,L)$ directly. 

\subsection{Approach}\label{approach}
In this article we prove Theorem \ref{power} first, and derive other Theorems from this theorem.
The main ingredients in our proof consist of the algebraic cobordism group $\SL$, the enumerative number $d_r(S,L)$ and the moduli stack of families of ideal sheaves. 
These techniques are developed in order to study the degeneration of line bundles on surfaces, and to derive a degeneration formula for the generating functions. 
As a consequence, we will prove when $L$ is $(5r-1)$-very ample, the number of nodal curves in $|L|$ on $S$ only depend on the class of $[S,L]$ in $\SL$, which can be computed by $L^2$, $LK$, $c_1(S)^2$ and $c_2(S)$. 

The algebraic cobordism theory has been developed by Levine and Pandharipande \cite{LP}. 
They call \begin{align*}[X_0]-[X_1]-[X_2]+[X_3]\end{align*} a \textit{double point relation} if there exists a flat family of projective schemes $\pi: \cX \to \PP^1$ satisfying the following properties: firstly, $\cX$ is smooth and $X_0$ is the fiber over $0\in \PP^1$ and  is smooth. Secondly, the fiber over $\infty\in \PP^1$  is the union of two smooth components $X_1$ and $X_2$ intersecting transversally along a smooth divisor $D$. Thirdly, $X_3= \PP(1_D\oplus N_{X_1/D})$ is a $\PP^1$ bundle over $D$. 
Define the addition of two schemes to be the disjoint union and the multiplication to be the Cartesian product. 
The algebraic cobordism ring is defined to be the ring generated by all smooth projective schemes modulo the double point relation. 
 
Since the problem of counting nodal curves is about a surface $S$ and a line bundle $L$ on $S$, we generalize Levine and Pandharipande's construction to pairs of line bundles on surfaces. Let $L_i$ be line bundles on $X_i$, we call 
\begin{align*}[X_0,L_0]-[X_1, L_1]-[X_2, L_2]+[X_3, L_3]\end{align*}
an \textit{extended double point relation} if $[X_0]-[X_1]-[X_2]+[X_3]$ is a double point relation, and there exists a line bundle $\cL$ on $\cX$ such that $L_i=\cL|_{X_i}$ for $i=0,1,2$; $L_3=\eta^* (\cL|_D)$ where $\eta : X_3\to D$ is the projection. 

In \cite{LP}, Levine and \Pan\ defined the \textit{algebraic cobordism group of surfaces and line bundles} $\SL$ to be the vector space over $\QQ$ spanned by all pairs $[S,L]$ modulo all extended double point relations.
The subscript $(2, 1)$ captures the dimension of surfaces and the rank of line bundles.
In Section \ref{chap3}, we prove $\SL$ is a four-dimensional vector space over $\QQ$ and $(L^2, LK, c_1(S)^2, c_2(S))$ induces the isomorphism from $\SL$ to $\QQ^4$. 
Consequently, the class of $[S,L]$ in $\SL$ is linear in $L^2, LK, c_1(S)^2$ and $c_2(S)$.  
Two important bases of $\SL$ are
$$\{[\PP^2, \cO], [\PP^2, \cO(1)], [\ptp, \cO], [\ptp, \cO(1,0)]\}  \text{ and }$$
$$\{ [\PP^2, \cO], [\PP^2, \cO(1)], [S_1, L_1],[S_2, L_2] \}, $$
where $S_i$ are K3 surfaces and $L_i$ are primitive classes on $S_i$ with $L_1^2\neq L_2^2$. 

In general, for pairs of vector bundles of rank $r$ on smooth projective scheme of dimension $n$, one can define the algebraic cobordism group $\omega_{n,r}$. Lee and \Pand\ \cite{LeeP} has studied the structure of $\omega_{n,r}$. Our result about $\SL$, although written independently, is a special case of \cite{LeeP}. 

Another important ingredient of our proof is the enumerative number $d_r(S,L)$. This number is defined by  G\"{o}ttsche (\cite{Gott}) is the intersection number of a closed subscheme $W^{3r}$ of $S^{[3r]}$ and the $2r^{th}$ Chern class of the tautological bundle of $L$ on $S^{[3r]}$ i.e.
$$d_r(S,L)=\int_{W^{3r}} c_{2r}(L^{[3r]}). $$ 
He proved that when $L$ is $(5r-1)$-very ample, $d_r(S,L)$ equals the number of  $r$-nodal curves in $[S, L]$. 
Although it is difficult to compute $d_r(S,L)$ directly, we believe it is the correct object to investigate for two reasons. 

First, when a pair of smooth surface and ample line bundle degenerates to the singular fiber, ampleness of the line bundle is usually not preserved. 
Without ampleness, the universal polynomial $T_r(L^2, LK, c_1(S)^2, c_2(S))$ do not necessarily equal to the number of $r$-nodal curves in $[S,L]$. 
On the contrary, the numbers $d_r(S,L)$ can be defined for every line bundle, even if it is not ample. 
In the end, we will show that for every pair $[S,L]$, $d_r(S,L)$ is the universal polynomial $T_r(\fourtop)$. 

Second, since $d_r(S,L)$ is defined as an intersection number, it has many good properties under degeneration. 
Suppose $\pi: \cX \to \PP^1$ defines a double point relation and $U$ is a Zariski open set of $\PP^1$ such that all fibers are smooth except $\pi^{-1}(\infty)$. 
J.~Li and B.~Wu \cite{Wu} constructed a moduli stack $\cXn\to U$, which is the moduli stack of rank one stable relative ideal sheaves in the family $\cX_U=\cX\times_{\PP^1} U \to U$. 
For each integer $n$, the moduli stack $\cXn \to U$ can be viewed as a family of Hilbert schemes of $n$ points on fibers. 
The generic fiber of $\cXn\to U$ at $\infty \neq t\in U$  is the Hilbert scheme $X_t^{[n]}$ and the special fiber is the union of products of relative Hilbert schemes
$$\bigcup_{k=0}^n (X_1/D)^{[k]}\times (X_2/D)^{[n-k]}. $$ 

For each $r\in \NN$, recall $d_r(S,L)$ is defined to be $\ds\int_{W^{3r}} c_{2r}(L^{[3r]})$ . The closed subscheme $W^{3r}$ in $S^{[3r]}$ can be extended globally to a family of closed subschemes $\cW^{3r}$ in $\cX^{[3r]}$.
In addition, if $\cL$ is a line bundle on $\cX$, then its restriction on $U$ similarly defines a tautological bundle $\cLn$ on $\cXn$ for each $n$. 
Thus on $\cX^{[3r]}$, the intersection of $c_{2r}(\cL^{[3r]})$ and $\cW^{[3r]}$ defines a family of zero cycles. Consider an extended double point relation $[X_0,L_0]-[X_1, L_1]-[X_2, L_2]+[X_3, L_3]$.  
By rational equivalence of the fibers of $W^{3r}$ over $0$ and $\infty$, we show the generating function
$$\phi(S,L)(x)=\sum_{r=0}^{\infty} d_r(S,L)x^r$$ satisfies the following degeneration formula
\begin{align*}\phi(X_0,L_0)=\phi(X_1/D,L_1)\phi(X_2/D,L_2)\end{align*}
where $\phi(X_i/D,L_i)$ is the ``relative'' generating function. 

The relation of absolute generating function $\phi$ can be computed by applying the degeneration formula on several families to eliminate the relative functions, which is
$$\phi(X_0,L_0)=\frac{\phi(X_1,L_1)\phi(X_2,L_2)}{\phi(X_3, L_3)}. $$ 
Thus $\phi$ induces a group homomorphism from the algebraic cobordism group $\SL$ to $(\qqx,\, \cdot\,)$.   
Consequently, Theorem \ref{power} is proved and $d_r(S,L)$ equals the universal polynomial $T_r(L^2, LK, c_1(S)^2, c_2(S))$ for all smooth projective surfaces $S$ and line bundles $L$. 
Since $d_r(S,L)$ equals the number of $r$-nodal curves when $L$ is $(5r-1)$-very ample, Theorem \ref{univpoly} is proved as a corollary. 

The generating function in Theorem \ref{DG2} is 
$$\gamma(S,L)(q)=\sum_{r\in \ZZ}T_r(L^2, LK, c_1(S)^2, c_2(S))(DG_2(\tau))^r= \phi(S,L)(DG_2).$$ 
Therefore this generating function also induces a homomorphism from $\SL$ to $\qqq$. 
Since Bryan and Leung \cite{BL} have computed $\gamma(S,L)(q)$ on generic K3 surfaces and primitive classes, we shall use a different basis 
$$\{[\PP^2, \cO], [\PP^2, \cO(1)], [S_1, L_1], [S_2, L_2]\}$$
of $\SL$,  where $S_i$ are K3 surfaces and $L_i$ are primitive classes on $S_i$ with $L_1^2\neq L_2^2$. 
By Theorem \ref{power}, $\gamma(S,L)(q)$ is a weighted product of $\gamma(\PP^2, \cO)$, $\gamma(\PP^2, \cO(1))$, $\gamma(S_1, L_1)$ and $\gamma(S_2, L_2)$. This proves Theorem \ref{DG2}. See Section \ref{chap4} for more details about proofs and computation. 

\subsection{Outline} 
In Section \ref{chap3} we construct the algebraic cobordism group $\SL$ and study its structure. 
We prove $\SL$ is a four dimensional vector space and $L^2$, $LK$, $c_1(S)^2$, $c_2(S)$ are invariants of extended double point relations. 
Therefore we can describe all bases of $\SL$ and degenerate any $[S,L]$ to basis elements. 

Section \ref{degformula} is dedicated to the enumerative number $d_r(S,L)$ and its generating function $\phi(S,L)(x)$. The main result in this section is a degeneration formula about $\phi(S,L)(x)$ for pairs satisfying an extended double point relation.   

Finally, in Section \ref{chap4} we combine the techniques developed in Section \ref{chap3} and \ref{degformula} to prove Theorem \ref{univpoly}, \ref{DG2} and Theorem \ref{power}. 
We express the infinite series $A_i$ and $B_i$ as the weighted product of generating functions on $\PP^2$, $\ptp$, and K3 surfaces.

\subsection{Acknowledgments} I am very grateful to my advisor Jun Li for his
supervision, advice, and guidance during my Ph.D. years in Stanford. His contribution of time, ideas and patience is invaluable to me. Special thanks also to Ravi Vakil for teaching me lots of things and great contribution to Stanford's algebraic geometry group. 
Moreover, I thank Lothar G\"{o}ttsche for the explanation of his approach and Rahul Pandharipande for useful suggestions. % and suggestion about Lemma \ref{new} (1)

\subsection{Notation and Convention}

\begin{enumerate}
  \item All surfaces are assumed to be complex, projective, smooth and algebraic. 
	\item An $r$-nodal curve is a reduced connected curve that has exactly $r$ nodes and no other types of singularity. 
		\item We always denote by $[S,L]$ a pair of a smooth projective surface $S$ and a line bundle $L$ on $S$. 
	\item The number of $r$-nodal curves in a pair $[S,L]$ means the number of $r$-nodal curves in a generic $r$-dimensional linear system of $|L|$ on $S$.

\end{enumerate}

\section{The algebraic cobordism group $\SL$}\label{chap3}

\subsection{Outline}
Degeneration methods have been applied widely in algebraic geometry. 
The general principle is that, if desired properties on general spaces is preserved under degenerations, 
then it can be verified by studying well-understood spaces.

To count the number of nodal curves, only degenerating surfaces is not sufficient.  
Instead, we consider the degeneration of pairs $[S,L]$ where $S$ is a smooth projective surface and $L$ is a line bundle on $S$. 
In Section 14.4 of \cite{LP},  Levine and Pandharipande defined the algebraic cobordism theory  $\SL$ for surfaces and line bundles. In this section we will find its dimension, bases and invariants\footnote{The algebraic cobordism theory developed in this section can be carried out for surfaces are over any field of characteristic zero.}. 
The study of $\SL$ has been extended to schemes of any dimension and vector bundles of arbitrary rank in \cite{LeeP}. Therefore the result in this section become a special case of \cite{LeeP}. We keep our proof in order to make this article self-contained and provide an elementary argument.

\subsection{Algebraic cobordism}
\begin{defn}\label{dp} Suppose $[X_i, L_i]$ are pairs of smooth projective surfaces and line bundles for $i=0,1,2,3$.  
The \textit{extended double point relation} is defined by
\begin{align}[X_0,L_0]=[X_1, L_1]+[X_2, L_2]-[X_3, L_3]\label{eq:dp} \end{align}with the assumption that 
there exists  a flat family of surfaces $\pi: \cX \to \PP^1$ and a line bundle $\cL$  on $\cX$  which satisfies the following properties: 
\begin{enumerate}
	\item $\pi^{-1}(\infty)= X_1 \cup_{D} X_2$  is a union of two irreducible smooth components, and they intersect transversally along a smooth divisor $D$;
	\item $\cX$ is smooth and of pure dimension three; the morphism $\pi$ is smooth away from the fiber $\pi^{-1}(\infty)$; 
	\item the fiber over $0\in \PP^1$ equals $X_0$, by (2) it is a irreducible smooth projective surface;
	\item $L_i$ is the restriction of $\cL$ on $X_i$ for  $i=0, 1, 2$;
	
	\item $X_3$ is $\PP(\cO_D\oplus N_{X_1/D})\cong \PP(N_{X_2/D}\oplus \cO_D)$, $\eta :X_3 \to D$ is the bundle morphism, and $L_3=\eta^* (\cL|_D)$. 
\end{enumerate}
 
\end{defn}

\begin{rem}We use $X_1/D$ and $X_2/D$ to emphasize the divisor $D$ when discussing  the relative geometry on $X_i$. \end{rem}

Let $M$ be the $\QQ$-vector space\footnote{We use $\QQ$ here because it is enough for our purpose and for simplicity. In general, the coefficients can be $\ZZ$ as discussed in \cite{LP}. } spanned by pairs of smooth projective surfaces and line bundles, and 
let $R$ be the subgroup spanned by all extended double point relations. 
We define the \textit{algebraic cobordism group of surfaces and line bundles} to be $$\SL=M/ R.$$ 
 
By definition $\SL$ is a $\QQ$-vector space. We will find its dimension, its bases, and the invariants of this degeneration theory. 
To find a basis of $\SL$, the central idea is degenerating the pairs of surfaces and line bundles $[S, L]$ into the sum of ``simpler'' pairs until the surfaces of generators become $\ptp$ and $\PP^2$. 
When Levine and Pardharipante \cite{LP} first defined this group $\SL$ (which they  called $\omega_{2,1}(\CC)$), 
They found an infinite set of generators which consisted of all line bundles on $\PP^2$, $\ptp$ and the blow-up of $\PP^2$ at a point.  
We will refine this result by showing $\SL$ is a four-dimensional vector space and it has a basis
$$\{[\PP^2, \cO], [\PP^2, \cO(1)], [\ptp, \cO], [\ptp, \cO(1,0)]\}.$$
During degeneration, many geometric properties are forgotten except four topological numbers $L^2, LK, c_1(S)^2$ and $c_2(S)$, and they are the only four invariants of $\SL$. 
As a result, the class of $[S,L]$ in $\SL$ only depends on $L^2$, $LK$, $c_1(S)^2$ and $c_2(S)$. 
A simple criterion to determine all bases of $\SL$ using these four numbers will be given in the end of this section.

First we review the algebraic cobordism theory. In \cite{LP}, Levine and \Pan\ consider projective morphisms 
$$\pi: Y\to X\times \PP^1$$ and composition $$\pi_2: p_2\circ \pi: Y\to \PP^1\,$$ where $k$ is a field of characteristic 0, $Y$ is a smooth quasi-projective scheme of pure dimension and $X$ is a separated scheme of finite type over $k$. 
Assuming $\pi_2^{-1}(\infty)=A\cup B$ where $A$ and $B$ are smooth Cartier divisors intersecting transversely along 
$D=A\cap B$, $$X_3=\PP(\cO_D\otimes N_{A/D})= \PP( N_{B/D}\otimes\cO_D)$$ is a $\PP^1$ bundle over $D$, and $Y_0$ is the fiber of $\pi_2$ over $0$. 
The double point relation over $X$ defined by $\pi$ is  
$$[Y_0\to X]-[A\to X]-[B\to X]+[X_3\to X]\,.$$

Let $\cM(X)^+$ be the free additive group generated by $[M \to X]$, where $M$ is a quasi-projective smooth scheme and the morphism is projective. 
Denote $\cR(X)$ by the subgroup generated by all double point relations over $X$ and $\omega_*(X)=\cM(X)^+/\cR(X)$. Then the following theorem states $\omega_*(X)$ is  isomorphic to the other algebraic cobordism ring $\Omega_*(X)$, which is constructed in \cite{LM}.

\begin{theo}[\textup Levine and Pandharipande \cite{LP}\textup]\label{LP}
There is a canonical isomorphism $$\omega_*(X)\cong \Omega_*(X)\,.$$\end{theo}      

In particular, if $k$ is a field, we denote $\Omega(k)$ to be $\Omega_*(\text{Spec }k)$ and  
$\omega_*(k)$ to be $\omega_*(\text{Spec }k)$ in short.
In \cite{LM}, Levine and Morel showed that  $\Omega_*(k)$ is isomorphic to the Lazard ring $\LL_*$; it follows that $\LL_*$ is also isomorphic to $\omega_*(k)$. 
Furthermore, it is well known that $\LL_* \otimes_{\ZZ} \QQ$ has a basis formed by products of projective spaces. 
Hence

\begin{coro}
$$\omega_*(k) \otimes_{\ZZ} \QQ=
\bigoplus_{\lambda=(\lambda_1, \ldots, \lambda_r)} \QQ[\PP^{\lambda_1}\times \ldots\times\PP^{\lambda_{r}}]$$ where the index $\lambda$ belongs to $\NN^r$ for some positive integer $r$.  
In particular, $\omega_2(k)$ is generated by $\PP^2$ and $\ptp$ over $\QQ$, i.e. all smooth projective surfaces can be degenerated to the sum of $\PP^2$ and $\ptp$ using algebraic cobordism. 
\end{coro}

We only consider the case that $X=\text{Spec }\CC$. It is easy to see that the extended double point relation generalizes the double point relation in dimension two.

\subsection{Theorem}
The following theorem is the main result of this section: 

\begin{theo}\label{slfl}Let $S$ be a smooth projective surface and $L$ be a line bundle on $S$, then 
\begin{align}\label{sl}[S,L]= a_1[\PP^2, \cO]+a_2[\PP^2, \cO(1)]+a_3 [\ptp, \cO]+a_4[\ptp, \cO(1,0)]  \end{align}
in $\SL$, where 
\begin{align*}
a_1&= -L^2+\frac{c_1(S)^2+c_2(S)}{3}-c_2(S), &a_2&= L^2,\\
a_3&=L^2+\frac{LK+L^2}{2}-\frac{c_1(S)^2+c_2(S)}{4}+c_2(S), & a_4&=-L^2-\frac{LK+L^2}{2}.
\end{align*}
In other words, $\SL$ is a four dimensional vector space spanned  by four elements 
$$[\PP^2, \cO], [\PP^2, \cO(1)], [\ptp, \cO] \text{ and }[\ptp, \cO(1,0)]\,.$$ Moreover, $(\fourtop)$ defines an isomorphism from $\SL$ to $\QQ^4$. 

\end{theo}

\begin{rem} By Theorem \ref{slfl}, the class of $[S,L]$ in $\SL$ is uniquely determined by topological numbers 
$( L^2, LK, c_1(S)^2, c_2(S))$. 
Noether's formula states $$\chi(\cO_S)=\frac{1}{12} (c_1(S)^2+c_2(S))$$ thus
$c_1(S)^2+c_2(S)$ is divisible by three and four. 
In addition, the Riemann-Roch formula for $L$ is $$\chi(L)=\chi(\cO_S)+\frac{1}{2}(L^2-LK)$$ thus $L^2+LK$ is divisible by two. Therefore all coefficients in \eqref{sl} are integers. \end{rem}

Before we prove Theorem \ref{slfl}, several lemmas are needed. 
\begin{notation}Suppose $C$ is a smooth curve and $N$ is a line bundle on $C$. We call $\PP_N$ to be the $\PP^1$ bundle $\PP_C(\cO_C\oplus N)$ over $C$ and write $h$ for the class of the tautological bundle $\cO_{\PP_N}(1)$ in $\mathrm{Pic}\, \PP_N$. 
Recall $\text{Pic}\, \PP_N =\eta^* \text{Pic}\,C\oplus \ZZ h$, where $\eta$ be the structure morphism $\PP_N \to C$. 
Moreover, denote the Hirzebruch surfaces $\PP_{\PP^1}(\cO_{\PP^1}\oplus \cO_{\PP^1}(k))$ by $\FF_k$.  
Then $\text{Pic } \FF_k \cong \ZZ h\oplus \ZZ f$, here $f$ is the fiber class and $h$ is the class of $\cO_{\FF_k}(1)$. 
Their intersection numbers are $$h^2=k,\ hf=1, \ f^2=0\,.$$

We will abuse notation slightly by using $\eta$ to denote all structure maps from projective bundle to its base curve. The referred morphism should be clear from the text.  \end{notation}

\begin{lemm}\label{4.1} Let $C$ be a smooth curve in $S$, $L$ be a line bundle on $S$, $N$ be the normal bundle of $C$ in S and $\eta: \PP_N\to C$ be the structure map. 
Then for any integer  $k$,  
\begin{align*}[S,L]=&[S, L\otimes \cO_S(-kC)]+[\PP_N, \eta^* \left(L|_C\right) \otimes \cO_{\PP_N}(kC)]
-[\PP_N, \eta^*((\cO_S(-kC)\otimes L)|_C)]\,\end{align*} in $\SL$. 
\end{lemm}

\begin{proof} Let $\mathrm{bl}: \cX\to S\times \PP^1$ be the blowup $S\times \PP^1$ along $C \times  \{\infty\}$, $\mathrm{pr_1}$ be the projection from $S\times \PP^1$ to $S$ and $\pi: \cX \to S\times \PP^1 \to \PP^1$  be the composition of $\mathrm{bl}$ and projection to $\PP^1$. Then $\cX$ is a family of surfaces with general fiber $S$ and singular fiber $\pi^{-1}(\infty)= S \cup_C \PP_N$. 
In addition, let $$\cL=\mathrm{bl}^*\mathrm{pr}_1^*L\otimes \cO_{\cX}(-kE)\,,$$ where $E\cong\PP_N$ is the exceptional divisor of the blow up morphism. 
Then the desired formula is the extended double point relation on $(\cX, \cL)$.  
\end{proof}

\begin{coro}\label{4.15} Every pair $[S,L]$ in $\SL$ is equal to the sum of $[S, \cO]$ and pairs on ruled surfaces of the form $\PP_N$. 
\end{coro}

\begin{proof}
Lemma \ref{4.1} implies that we can twist the line bundle $L$ by effective divisors so that the remaining terms are pairs on $\PP_N$. Since every line bundle on a projective surface can be written as the difference of two very ample line bundles, the statement can be proved by applying Lemma \ref{4.1} with $k=1$ twice . 
\end{proof}

In Levine and Pandharipande's algebraic cobordism theory, the classes of surfaces are generated freely by $\PP^2$ and $\ptp$ over $\QQ$. Moreover, the class of $S$ depends only on $c_1(S)^2$ and $c_2(S)$ and is given by  
$$[S]=\left(\frac{c_1(S)^2+c_2(S)}{3}-c_2(S)\right) [\PP^2]
+ \left(-\frac{c_1(S)^2+c_2(S)}{4}+c_2(S)\right) [\ptp] \,.$$ 
Since extended double point relations with trivial line bundles reduce to double point relations, 
similar formula also holds in $\SL$:

\begin{lemm}\label{4.2} $$[S,\cO]=\left(\frac{c_1(S)^2+c_2(S)}{3}-c_2(S)\right) [\PP^2, \cO]
+ \left(-\frac{c_1(S)^2+c_2(S)}{4}+c_2(S)\right) [\ptp, \cO]  $$ in $\SL$. \end{lemm}

According to Lemma \ref{4.1} and \ref{4.2}, $\SL$ is generated by pairs on $\PP^2$ and ruled surfaces of the form $\PP_N$. 
Consequently, in the following lemmas we concentrate on degenerating these two kinds of pairs. A useful observation is 
that deformation of ruled surfaces can be achieved by deforming the base curve $C$ and constructing a family of $\PP^1$-bundle on each fiber. 
 
\begin{lemm}\label{PClemma} Suppose $N$ and  $N'$ are two line  bundles over a smooth curve $C$ and $x$ is a closed point on C, $n \in \NN$. 
Then
\begin{align*}[\PP_{N'}, \eta^*N\otimes \cO_{\PP_{N'}}(n)]
=& [\PP_{N'(-x)}, \eta^*N\otimes \cO_{\PP_{N'(-x)}}(n)] + [\FF_1, \cO_{\FF_1}(n)]- [\ptp, \cO(n,0) ]\end{align*} in $\SL$. 
\end{lemm}

\begin{proof}
Define $X$ to be the blow-up of $C\times \PP^1 $ at the closed point $x\times \infty$. 
$X$ is a family of curves whose general fiber is $C$ and whose fiber over $\infty$ is isomorphic to the union of $C$ and exceptional curve $E\cong \PP^1$. $C$ and $E$ intersect at the point $x$. 
There are morphisms
$$X \stackrel{\mathrm{bl}}{\ra} C\times \PP^1  \stackrel{\mathrm{pr_1}}{\ra} C,$$
where $\mathrm{bl}: X\to  C\times \PP^1$ is the blow-up morphism and $\mathrm{pr_1}:  C\times \PP^1 \to C$ is the projection. 

Let $M'$ be the pullback line bundle $\mathrm{bl}^*\mathrm{pr}_1^*N'$ and consider the line bundle $M'(-E)$ on $X$. 
The restriction of $M'(-E)$ on general fiber is $N'$, on the first component of singular fiber $C$ is $N'(-x)$, on $\PP^1$ is $\cO(1)$, and it restricts to trivial line bundle at $x$. 

Therefore, the $\PP^1$-bundle $\cX=\PP_X(\cO\oplus M'(-E))$ over  $X$ is a family of surfaces with generic fiber 
$\PP_C(\cO \oplus N' )$ and  singular fiber  $\PP_C(\cO \oplus N(-x))\cup_{\PP^1}\FF_1$. 

Let $\pi: \cX \to X$ be the bundle map. The morphism
$$\cX\stackrel{\pi}{\ra} X \stackrel{\mathrm{bl}}{\ra} C\times \PP^1 \stackrel{\mathrm{pr_2}}{\ra}\PP^1$$ 
defines $\cX$ as a flat family of surfaces over $\PP^1$.
Denote the generic fiber $\PP_C(\cO \oplus N')$ by $X_0$, components of the singular fiber $\PP_C(\cO \oplus  N'(-x))$ by $X_1$, $\FF_1$ by $X_2$ and their intersection $\PP^1$ by $D$. 
Then $\cX$ is a family of surfaces which satisfies the conditions in Definition \ref{dp}. 

Next we construct a line bundle $\cL$ on $\cX$. Let $\cO_{\cX}(1)$ be the tautological bundle of the projective bundle 
$\cX$ and  $\cL$ be $\pi^*\mathrm{bl}^*\mathrm{pr}_1^* N \otimes \cO_{\cX}(n)$. The restriction of $\cL$ on each fiber is:
$$\cL|_{X_0} \cong \eta^*N\otimes \cO_{\PP_{N'}}(n), \ \cL|_{X_1} =\eta^*N\otimes \cO_{\PP_{N'(-x)}}(n),\ 
\cL|_{X_2}=\cO_{\FF_1}(n),\  \cL|_D= \cO_{\PP^1}(n).$$  
Then the  desired formula is the extended double point relation on $\cX$ and $\cL$.  
\end{proof}

Every line bundle on $\PP_{N'}$ can be written as $\eta^* N\otimes \cO_{\PP_{N'}}(n)$.  
We can write $N'$ as $\cO(\sum a_ix_i-b_iy_i)$ for some closed points $x_i$, $y_i$ on $C$ and positive integers $a_i$, $b_i$. By using Lemma \ref{PClemma}  inductively on $\sum|a_i|+\sum|b_i|$, pairs on $\PP_{N'}$ can be reduced to pairs on $\PP_C(\cO_C\oplus \cO_C)\cong C\times \PP^1$ with excess terms on $\FF_1$ and $\ptp$.
 
Furthermore, line bundles on $C\times \PP^1$ can also be reduced using the following lemma:

\begin{lemm}\label{CP1} Every pair on $C\times \PP^1$ can be written as the sum of pairs on $\ptp$. 
\end{lemm} 
\begin{proof} Every line bundle on $C\times \PP^1$ is isomorphic to $L\boxtimes \cO(n)$ for some integer $n$ and line bundle $L$ on $C$. 
Let $x$ be a closed point on $C$ and  $\cX$ be the blow-up of $C\times \ptp$ at $x\times \PP^1\times  \infty $.
The third projection $\cX\to \PP^1$ defines $\cX$ as a family of surface over $\PP^1$ with general fiber 
$X_0\cong C\times\PP^1$ and singular fiber $X_1\cup_D X_2=(C \times \PP^1) \cup_{\PP^1}(\ptp)$. 

Let $\cL$ be $\mathrm{pr}_1^*L \otimes \mathrm{pr}_2^*\cO_{\PP^1}(n)\otimes \cO_{\cX}(-X_2)$, where $\mathrm{pr_i}$ are projections to $C$ and the first $\PP^1$. 

The extended double point relation from $\cX$ and $\cL$ is 
\begin{equation}\begin{split}\label{eq:CP1}[C\times \PP^1, L\boxtimes\cO(n)]&=[C\times \PP^1, L(-x)\boxtimes\cO(n)]+[\ptp, \cO(1,n)]\\
&-[\ptp, \cO(0,n)]\end{split}\end{equation}

By applying \eqref{eq:CP1} several times, we can see that $[C\times \PP^1, L\boxtimes\cO(n)]$ is the sum of pairs on $\ptp$ and  $[C\times \PP^1, \cO_C\boxtimes\cO(n)]$. 

On the other hand, if $$[C_0]=[C_1]+[C_2]-[C_3]$$ is a double point relation of curves defined by a family $Y$, then ($Y\times \PP^1, \pr_2^* \cO_{\PP^1}(n))$ gives an extended double point relation
$$[C_0\times \PP^1, \pr_2^* \cO_{\PP^1}(n)]=[C_1\times \PP^1, \pr_2^* \cO_{\PP^1}(n)]+[C_2\times \PP^1, \pr_2^* \cO_{\PP^1}(n)]-[C_3\times \PP^1, \pr_2^* \cO_{\PP^1}(n)].$$
Since every smooth curve can be degenerated to several $\PP^1$'s using double point relation, 
$[C\times \PP^1, L\boxtimes\cO(n)]$ is the sum of pairs on $\ptp$. 

\end{proof}

From the previous lemmas, we have shown that $\SL$ is spanned by $[\PP^2, \cO]$, pairs on $\ptp$ and 
$[\FF_1, \cO(nh)]$. 

\begin{lemm}\label{new} The following equalities hold in $\SL$ for all integers $a$, $b$. 
\begin{enumerate}
   \item $[\ptp, \cO(a,a)]=2[\FF_1, ah]-[\ptp, \cO(0,a)]$  \smallskip
   \item $[\ptp, \cO(a,b)]=ab[\ptp, \cO(1,1)]-(2ab-a-b)[\ptp, \cO(0,1)]\\
   \phantom{[\ptp, \cO(a,b)=}+(a-1)(b-1)[\ptp, \cO]$\smallskip
   \item $[\PP^2, \cO(1)]=[\FF_1, h]+[\PP^2, \cO]-[\ptp, \cO]$ \smallskip\end{enumerate}\end{lemm}

\begin{proof}
\begin{enumerate}
\item Fix two points $0, \infty \in \PP^1$, consider the family $\mathrm{bl}_{\, 0\times \infty} \ptp \stackrel{\mathrm{pr}_2}{\longrightarrow} \PP^1$. 
The fiber over $0$ is $C_0\cong \PP^1$; fiber over $\infty$ is the union of $C_1$ and $C_2$, both isomorphic to  $\PP^1$ and they intersect at one point. A section of the family either passes through $C_1$ or $C_2$, hence induce two line bundles $M_1$ and $M_2$ on $\mathrm{bl}_{\, 0\times \infty} \ptp$ which satisfy: 
\begin{align*}
&M_1|_{C_0}= \cO(1),\ &M_1|_{C_1}&= \cO(1),\  &M_1|_{C_2}&= \cO;&\\
&M_2|_{C_0}= \cO(1),\ &M_2|_{C_1}&= \cO,   \  &M_2|_{C_2}&= \cO(1).&\\
\end{align*}   
Then the composition $$\PP(M_1 \oplus M_2)\to \mathrm{bl}_{\, 0\times \infty} \ptp \stackrel{\mathrm{pr}_2}{\longrightarrow} \PP^1$$  defines $\PP(M_1 \oplus M_2)$ as a family of surfaces. 
Since  $\PP(M_1 \oplus M_2)$ is a $\PP^1$-bundle over $\mathrm{bl}_{\, 0\times \infty} \ptp$, the tautological bundle $\cO_{\PP(M_1 \oplus M_2)}(1)$ is naturally defined. 
The extended double point relation from $\PP(M_1 \oplus M_2)$ and $\cO_{\PP(M_1 \oplus M_2)}(a)$ is exactly 
$$[\ptp, \cO(a,a)]=2[\FF_1, ah]-[\ptp, \cO(0,a)]. $$

\item 
Apply \eqref{eq:CP1} to the case $C=\PP^1$, $L=\cO(a)$ and use induction on $|a|$, we have
$$[\ptp, \cO(a,b)]=a[\ptp, \cO(1,b)]-(a-1)[\ptp, \cO(0,b)].$$ 
By using the identity on both $a$ and $b$, we get
\begin{align*}&[\ptp, \cO(a,b)]=a[\ptp, \cO(1,b)]-(a-1)[\ptp, \cO(0,b)]\\
 =& ab[\ptp, \cO(1,1)]-(2ab-a-b)[\ptp, \cO(0,1)]+(a-1)(b-1)[\ptp, \cO]. \end{align*}
\item Let $\cX$ be the blowup of $\PP^2 \times \PP^1$ along $pt \times \infty$ and $\cL$ be the the pullback of $\cO_{\PP^2}(1)$ in $\cX \stackrel{\bl}{\rightarrow}  \PP^2 \times \PP^1 \stackrel{\pr_2}{\rightarrow} \PP^2$. 
The identity we need to prove is the extended double point relation on $(\cX, \cL)$ plus 
\begin{enumerate}
	\item $\FF_1$ is isomorphic to $\PP^2$ with a point blown up;
	\item  because the Chern numbers of $\FF_1$ and $\ptp$ are the same, $[\FF_1]= [\ptp]$ in the algebraic cobordism ring of schemes, hence $[\FF_1, \cO]= [\ptp, \cO]$ in $\SL$. 
\end{enumerate}
\end{enumerate}\end{proof}

\begin{coro}\label{newcoro}$$[\ptp, \cO(1,1)]=2[\PP^2, \cO(1)]-2[\PP^2, \cO]+2[\ptp, \cO]-[\ptp, \cO(1,0)].$$\end{coro}
\begin{proof}This follows from Lemma \ref{new} (1) and (3). \end{proof}

Next, we show that $\SL$ has four independent invariants. 

\begin{prop}\label{slinv}
Suppose $$[X_0,L_0]=[X_1, L_1]+[X_2, L_2]-[X_3, L_3]$$ is an extended double point relation, then 
\begin{align*}
\left(L_0^2, L_0K_{X_0}, c_1(X_0)^2, c_2(X_0)\right)=\left(L_1^2, L_1K_{X_1}, c_1(X_1)^2, c_2(X_1)\right)\\
+\left(L_2^2, L_2K_{X_2}, c_1(X_2)^2, c_2(X_2)\right)
-\left(L_3^2, L_3K_{X_3}, c_1(X_3)^2, c_2(X_3)\right)\,.\end{align*} In other words, if 
$$\sum_i [X_i, L_i]=\sum_j [X_j, L_j] \ \  \text{ in } \SL$$ then 
$$\sum_i (L_i^2, L_iK_{X_i}, c_1(X_i)^2, c_2(X_i) )= \sum_j (L_j^2, L_jK_{X_j}, c_1(X_j)^2, c_2(X_j) ). $$
Thus $(\fourtop)$ induces an homomorphism from $\SL$ to $\QQ^4$.
\end{prop}
\begin{proof}
 
Recall that $X_3=\PP(\cO_D\oplus N_{X_1/D})$ and $L_3=\eta^* (\cL|_D)$. Simple calculation shows
\begin{align*} L_3^2&=0, \\
 L_3K_{X_3}&=(\deg(\cL|_D)f)(-2h+(\deg(N_{X_1/D})+2g(D)-2)f)=-2\deg(\cL|_D), \\
 K_{X_3}^2&=(-2h+(\deg(N_{X_1/D})+2g(D)-2)f)^2\\&=4 \deg(N_{X_1/D})-4  \deg(N_{X_1/D}) -8g(D)+8
 =-8g(D)+8.
\end{align*}
Because $X_0$ and $X_1\cup X_2$ are two fibers of $\cX\to \PP^1$, in the Chow group of $\cX$ 
$$[X_0]\cdot[X_0]=0, \ [X_1]\cdot[X_1]=-[X_2][X_1]=-[D]=[X_2]\cdot[X_2],$$ 
$$[X_0][X_1]=[X_0][X_2]=0\,.$$ 

Then 
\begin{align*}
L_0^2&=c_1(\cL)^2[X_0]= c_1(\cL)^2[X_1\cup_D X_2]= c_1(\cL)^2 [X_1]+c_1(\cL)^2[X_2]
=L_1^2+L_2^2-L_3^2,	\\
L_0K_{X_0}&= c_1(\cL)(c_1(K_{\cX})+[X_0])\cdot [X_0]= c_1(\cL)c_1(K_{\cX}) [X_0]  \\
&= c_1(\cL)c_1(K_{\cX}) [X_1\cup_D X_2]=c_1(\cL)c_1(K_{\cX})[X_1]+ c_1(\cL)c_1(K_{\cX})[X_2]\\
&=c_1(\cL)(c_1(K_{\cX})+[X_1]-[X_1]) [X_1]+ c_1(\cL)(c_1(K_{\cX})+[X_2]-[X_2]) [X_2]\\
&=L_1K_{X_1} -c_1(L_1)[X_1]^2 +L_2K_{X_2} -c_1(L_2)[X_2]^2\\
&= L_1K_{X_1} +c_1(L_1)[D] +L_2K_{X_2} +c_1(L_2) [D]\\
&=L_1K_{X_1} +L_2K_{X_2} +2\deg(\cL|_D)\\
&=L_1K_{X_1} +L_2K_{X_2} -L_3K_{X_3},\\
%%%%%%%%%%%%%%%%%%%%%%%%%
K_{X_0}^2&= (c_1(K_{\cX})+[X_0])^2\cdot[X_0]= c_1(K_{\cX})^2[X_1\cup_D X_2] \\
&= c_1(K_{\cX})^2 [X_1]+ c_1(K_{\cX})^2 [X_2]\\
&= (c_1(K_{\cX})+[X_1]-[X_1])^2\cdot[X_1]+  (c_1(K_{\cX})+[X_2]-[X_2])^2\cdot [X_2]\\
&=K_{X_1}^2+ K_{X_2}^2 -2(c_1(K_{\cX})+[X_1])\cdot[X_1]^2-2(c_1(K_{\cX})+[X_2])\cdot[X_2]^2+0\\
&=K_{X_1}^2+ K_{X_2}^2+2c_1(K_{X_1})[D]+ 2c_1(K_{X_2})[D]\\
&= K_{X_1}^2+ K_{X_2}^2 +4(2g(D)-2)\\
&=  K_{X_1}^2+ K_{X_2}^2- K_{X_3}^2.
\end{align*}
In addition,  the exact sequence
$$0\ra \cO_{X_1\cup_D X_2}\ra \cO_{X_1}\oplus \cO_{X_2} \ra \cO_{D} \ra 0$$
implies
\begin{align*}\chi(\cO_{X_0})&=\chi(\cO_{X_1\cup_D X_2})=\chi(\cO_{X_1})+\chi(\cO_{X_2})- \chi(\cO_{D})\\
&=\chi(\cO_{X_1})+\chi(\cO_{X_2})- (1-g(D)) =\chi(\cO_{X_1})+\chi(\cO_{X_2})-\chi(\cO_{X_3}). 
\end{align*} 
Using $\chi(\cO_S)=\frac{1}{12}(K_S^2+c_2(S))$,  we conclude
$$c_2(X_0)=c_2(X_1)+c_2(X_2)-c_2(X_3).$$
This completes the proof. 
\end{proof}

\begin{proof}[Proof of Theorem \ref{slfl}]
Although the lemmas above have implied that $\SL$ is spanned by $[\PP^2, \cO]$, $[\PP^2, \cO(1)]$, $[\ptp, \cO]$  and $[\ptp, \cO(1,0)]$, we outline our proof here as a summary. Here we don't keep track of the coefficients and use $*$ to indicate them. Recall that if $N$ is a line bundle on a smooth curve $C$, then $\PP_N$ is defined as the $\PP^1$ bundle $\PP_C(\cO_C\oplus N)$ on $C$.  

For every pair $[S,L]$ in $\SL$, 
\begin{align*}
&[S,L]\\= &[S,O]+ \text{ pairs on ruled surfaces } \PP_N  &\text{ (Corollary \ref{4.15})}\\
= &*[\PP^2, \cO]+ *[\ptp, \cO]+ \text{ pairs on ruled surfaces } \PP_N  &\text{ (Lemma \ref{4.2}) }\\
= & *[\PP^2, \cO]+ *[\ptp, \cO]+ \text{ pairs on } C\times \PP^1 + *\sum_{n\in \ZZ}\,[\FF_1, \cO_{\FF_1}(n)]&\text{ (Lemma \ref{PClemma})}\\
= & *[\PP^2, \cO]+ *[\ptp, \cO]+ \text{ pairs on } \ptp +*\sum_{n\in \ZZ}\,[\FF_1, \cO_{\FF_1}(n)] &\text{ (Lemma \ref{CP1})}\\
= &*[\PP^2, \cO]+ *[\ptp, \cO]+ *[\ptp, \cO(a,b)] &\text{ (Lemma \ref{new}(1))}\\
= & *[\PP^2, \cO]+ *[\ptp, \cO]+ *[\ptp, \cO(1,1)] + *[\ptp, \cO(0,1)] &\text{ (Lemma \ref{new}(2))}\\
= & *[\PP^2, \cO]+ *[\ptp, \cO]+ *[\PP^2, \cO(1)]+ *[\ptp, \cO(1,0)] &\text{ (Corollary \ref{newcoro})}
\end{align*}
As a result, $\SL$ are spanned by these four pairs. They are also independent because $\SL$ has four independent invariants \fourtopand (Proposition \ref{slinv}). Thus 
$$[\PP^2, \cO], [\PP^2, \cO(1)], [\ptp, \cO], [\ptp, \cO(0,1)]$$
is a basis of $\SL$. 

To prove equation \eqref{sl}, it suffices to find the coefficients. 
Suppose the class of $[S,L]$ in $\SL$ is
\begin{align*}[S,L]=a_1 [\PP^2, \cO]+a_2[\PP^2, \cO(1)]+ a_3[\ptp, \cO] +a_4[\ptp, \cO(1,0)]\,.\end{align*}
Since the topological numbers $(L^2, LK, c_1(S)^2, c_2(S))$ of 
$$[\PP^2, \cO], \,[\PP^2, \cO(1)], \,[\ptp, \cO] \text{ and } [\ptp, \cO(0,1)]\text{ are }$$ 
$$(0,0,9,3), \,(1, -3, 9,3), \,(0,0,8,4), \,(0,-2, 8,4)$$ and  $(L^2, LK, c_1(S)^2, c_2(S))$ are invariants of $\SL$, we have 
$$(L^2, LK, c_1(S)^2, c_2(S))=a_1 (0,0,9,3)+a_2(1, -3, 9,3)+a_3 (0,0,8,4)+a_4 (0,-2, 8,4). $$
Computation shows 
\begin{align*}
a_1&= -L^2+\frac{c_1(S)^2+c_2(S)}{3}-c_2(S), &a_2&= L^2,\\
a_3&=L^2+\frac{LK+L^2}{2}-\frac{c_1(S)^2+c_2(S)}{4}+c_2(S), & a_4&=-L^2-\frac{LK+L^2}{2}.\end{align*}
\end{proof}

\begin{rem} 
\begin{enumerate}
	\item From the proof of Theorem \ref{slfl}, we know a set of four elements 
	$$\{[S_i, L_i] \ | \ i=1,\ldots, 4\}$$ 
is a basis of $\SL$ if and only if the four vectors 
$$\{(L_i^2, L_iK_{S_i}, c_1(S_i)^2, c_2(S_i))\ |\ i=1, \ldots , 4\}$$ are linearly independent over $\RR$ (in fact over $\QQ$ is enough but we don't need it.)  
  \item Because the two vectors $(L^2, LK, c_1(S)^2, c_2(S))$ and $(LK, \chi(L), \chi(\cO), K^2)$ determine each other,   $(LK, \chi(L), \chi(\cO), K^2)$ is also an invariant of $\SL$. 
  As a result, if for a set of four elements in $\SL$, the corresponding vectors $(LK, \chi(L), \chi(\cO), K^2)$ are linearly independent, then the set is a basis of $\SL$. 
  \item  If $S_1, S_2$ are two K3 surfaces and $L_i$ are primitive classes on $S_i$ respectively. The four numbers  $(LK, \chi(L), \chi(\cO), K^2)$ of 
  $$\fB=\{ [\PP^2, \cO], [\PP^2, \cO(1)], [S_1, L_1], [S_2, L_2] \} \text{  are  }$$
  $$(0, 1, 1, 9), (-3, 3, 1, 9), \left(0, 2+\frac{L_1^2}{2}, 2, 0\right), \left(0, 2+\frac{L_2^2}{2}, 2, 0\right).$$ 
  It follows that $\fB$ is a basis if and only if $L_1^2\neq L_2^2$. 
\end{enumerate}
\end{rem}

\section{Degeneration Formula}\label{degformula}

\subsection{The enumerative number $d_r(S,L)$}
In this chapter we use the enumerative number $d_r(S,L)$ to study the number of nodal curves. 
This number $d_r(S,L)$ was defined  by \Gott\ \cite{Gott} for pairs of smooth projective surface and line bundle $[S, L]$ and he proved that $d_r(S,L)$ equals the number of $r$-nodal curves on $[S,L]$ if $L$ is $(5r-1)$-very ample. 
The goal of this chapter is to derive a degeneration formula for the generating function 
$$ \phi(S,L)(x)=\sum_{r=0}^{\infty} d_r(S,L)x^r$$ for pairs satisfying an extended double point relation.

Let $S^{[n]}$ be the Hilbert scheme of $n$ points on $S$, and let 
$Z_n\subset S\times S^{[n]}$ be the universal closed subscheme with projections 
$$p_n : Z_n \to S,\  q_n : Z_n \to S^{[n]}.$$ 
Define $L^{[n]} = {(q_n)}_* (p_n)^*L$. Because $q_n$ is finite and flat, $L^{[n]}$ is a vector bundle of rank $n$ on $S^{[n]}$. G\"{o}ttsche \cite{Gott} suggested the following approach using intersection numbers on Hilbert schemes: 

\begin{defn}[\cite{Gott} Definition 5.1] \label{S2r}
Let $W_0^{3r}$ be the locally closed subset 
$$\displaystyle \left\{\left.\coprod_{i=1}^{r}\text{Spec}(\cO_{S, x_i}/m^2_{S, x_i})\right|\  x_i \text{ are distinct closed points on }  S \right\}$$
and $W^{3r}\subset S^{[3r]}$ be the closure  of  $W_0^{3r}$ (with the reduced induced structure). 
It is easy to see that $W^{3r}$ is birational to $S^{[r]}$. Define $$d_r(S,L)=\int_{W^{3r}} c_{2r}(L^{[3r]}).$$ 
\end{defn}\smallskip

For simplicity, define $d_0(S,L)=1$ because the number of $0$-nodal curves (smooth) in a linear system is one.  

\begin{defn}We call $L$  \textsl{$k$-very ample} if for every zero-dimensional subscheme $\xi \subset S$ of length $k+1$, the natural map
$H^0(S,L)\to H^0(\xi, L \otimes \cO_{\xi})$ is surjective. \end{defn}
If $L$ and $M$ are very ample then $L^{\otimes k}\otimes M^{\otimes l}$ is $(k+l)$-very ample. 
In particular, very ampleness implies 1-very ampleness. 

We quote a result of G\"{o}ttsche below: 
\begin{prop}[\cite{Gott}, Proposition 5.2]\label{Gott'sthm}
Assume $S$ is a smooth algebraic surface and $L$ is a $(5r-1)$-very ample line bundle on $S$,
then a general $r$-dimensional sublinear system $V\subset |L|$ contains $d_r(S,L)$ curves with precisely $r$-nodes as singularities. 
\end{prop}

\subsection{Degeneration Formula}
In Section \ref{approach} we explained why $d_r(S,L)$ behaves better in families and is the right object to derive a degeneration formula. 
The goal of this section is to prove the degeneration formula for $d_r(S,L)$.  
More explicitly, we will show that if $$[X_0, L_0]=[X_1, L_1]+[X_2, L_2]-[X_3, L_3]$$ is an extended double point relation, then then number of nodal curves on $[X_0, L_0]$ can be determined by the numbers on $[X_1, L_1]$, $[X_2, L_2]$ and $[X_3, L_3]$. 
As a result, it is necessary to treat the number of curves with varied number of nodes together by considering the generating function
$$\phi(S,L)(x)= \sum_{r=0}^{\infty} d_r(S,L) x^r .$$

The following theorem is the main result of this section: 
\begin{theo}\label{dgfl} Suppose $[X_0,L_0]=[X_1, L_1]+[X_2, L_2]-[X_3, L_3]$ is an extended double point relation. Then  
$$\phi(X_0, L_0)=\frac{\phi(X_1, L_1)\cdot \phi(X_2, L_2)}{\phi(X_3, L_3)}. $$ 
In other words, $\phi$ is a homomorphism from $\SL$ to $(\qqx, \ \cdot \ )$. 
\end{theo}\smallskip

\subsection{Moduli stack of relative ideal sheaves}
The key tool is Jun Li and Baosen Wu's  \cite{Wu} construction of the moduli stack of stable relative ideal sheaves. Let $\infty \in C$ be a specialized point and  $\pi: X \to C$ be a flat projective family of schemes that satisfies
\begin{enumerate}
	\item $X$ is smooth and $\pi$ is smooth away from the fiber $\pi^{-1}(\infty)$; 
	\item $\pi^{-1}(\infty)=: X_1 \cup_{D} X_2$  is a union of two irreducible smooth components $X_1 $ and $X_2$ which intersect transversally along a smooth divisor $D$.
\end{enumerate}
In \cite{Wu}, Li and Wu defined the notion of a family of stable perfect ideal sheaves over $C$ and constructed $\fI_{X/C}^{\Gamma}$,  the moduli space of stable perfect ideal sheaves of type $\Gamma$ of $X\to C$. 
To make $\fI_{X/C}^{\Gamma}$  a stack, one has to replace $X$ by new spaces $X[n]$ so that $X$ and $X[n]$ have the same smooth fiber $X_t$ when $t\neq \infty$. Over $\infty$, the fiber of $X[n]$ is a semistable model 
$$X[n]_0=X_1 \cup \Delta_1 \cup \Delta_2\cup \ldots \Delta_{n-1}\cup X_2,$$ 
where $\Delta_i\cong \PP_D(\cO_D\oplus N_{X_1/D})$. 
The objects of $\fI_{X/C}^{\Gamma}$ are $(\cX/S, \cI)$ which consists of a family $\cX$ over a $C$-scheme $S$ and a family of stable ideal sheaves $\cI$ of type $\Gamma$ on $\cX/S$. The fibers $\cX_s$ of the family $\cX/S$ are required to be either smooth a fiber of $X/C$ or a semistable model $X[n]_0$ for some $n$. Under these settings, Li and Wu proved that the moduli space has many good properties:

\begin{theo}[\cite{Wu}]
The moduli stack $\fI_{X/C}^{\Gamma}$ is a separated and proper Deligne-Mumford stack of finite type over $C$. 
\end{theo} 
 
In our case, $\pi: X\to C$ is a family of surfaces and $\cI$ is a family of ideal sheaves of zero-dimensional closed subschemes of length $n$ . 
For $s\in S$, when $\cX_s$ is a smooth fiber of $X/C$, $\cI_s$ is automatically perfect and stable. If $\cX_s$ is a semistable model $X[n]_0$, then the support of $\cI_s$ can not lie on the intersection of two components of $X[n]_0$ and every component contains at least one point of the zero-set of $\cI_s$.

The fibers of $\pi^{[n]}: \fI_{X/C}^{n} \to C$ can be described as follows:
\begin{enumerate}
	\item The fiber over $\infty$  is the union of products 
$$\cup_{k=0}^n (X_1/D)^{[k]} \times (X_2/D)^{[n-k]}$$ for all possible $n\geq k\geq 0$. 
$(X_i/D)^{[n_i]}$ are the moduli spaces of stable relative ideal sheaves of $n_i$ points on $X_i/D$. $(X_i/D)^{[0]}=pt$. They are also  
separated and proper Deligne-Mumford stacks (see \cite{Wu}, Theorem 3.7). 
We denote this fiber  by $(X_1 \cup_{D} X_2)^{[n]}$. 
	\item $\fI_{X/C}^{n}$ is smooth and $\pi^{[n]}$ is smooth away from the fiber over $\infty$. 
	\item When $t\neq \infty$, the smooth fiber of ${\fI_{X/C}^{n}}$ over $t$ equals $X_t^{[n]}$ , the Hilbert schemes of $n$ points on $X_t$. 
\end{enumerate}

Therefore, $\fI_{X/C}^{n}$ can also be viewed as a family of Hilbert schemes of $n$ points on $X/C$. 
The following Definition-Proposition contains some facts proved in \cite{Wu}.

\begin{depr}\label{longdef}Let $\pi: \cX\to \PP^1$ be a family of surfaces described in Definition $\ref{dp}$ and $U$ is a Zariski open set of $\PP^1$ obtained by deleting those points with singular fibers except $\infty$, i.e. set-theoretically, 
$$U=\left\{t\in \PP^1\,|\, \text{ the fiber } \cX|_t \text{ is smooth }\right\} \cup \{\infty\}.$$ Then $\pi_U: \cX_U:=\cX\times_{\PP^1} U \to U$ is a family of surfaces with only one singular fiber $X_{\infty}$. 

\begin{enumerate}
	\item Define $\cX^{[n]}$ to be ${\fI_{{\cX_U}/U}^{n}}$ and $\pi^{[n]}: \cX^{[n]} \to U$ to be the structure morphism. 
	\item The fibers of $\pi^{[n]}$ are: 
	$${(\pin)}^{-1}(0)=X_0^{[n]} \text{ and } {(\pin)}^{-1}(\infty)= \cup_{k=0}^n (X_1/D)^{[k]} \times (X_2/D)^{[n-k]}.$$ The inclusions are denoted by
	$$i_0^{[n]}: X_0^{[n]} \to \cXn, \ \ i^{[k,n]}: (X_1/D)^{[k]} \times (X_2/D)^{[n-k]} \inclu \cXn.$$
	\item Let $Z_0^{[n]}$  be the universal closed subscheme  in $X_0\times X_0^{[n]}$. 
	For $i=1,2$, the universal closed subscheme of stable relative ideal sheaves on $X_i$ can also be constructed as a closed subscheme of  $(X_i/D)\times (X_i/D)^{[k]}$. We denote it by $Z_i^{[k]}$.
	
	\item There is a universal closed subscheme $\cZ\subset \cX_U\times \cX^{[n]}$ with structure morphisms 
	$$\xymatrix{
  \cZ  \ar[d]^{P} \ar[r]^{Q} & \cXn   \\
  \cX_U &  }  $$ 
  such that the composition $\pi_U \circ P:\cZ \to \cX_U \to U$ has fibers
	$$\cZ|_{0}=Z_0^{[n]}  \text{ and }$$
	$$\cZ|_{\infty}= \left(\cup_{k=0}^n Z_1^{[k]}\times (X_2/D)^{[n-k]}\right)\cup 
	\left( \cup_{k=0}^n  (X_1/D)^{[k]}\times Z_2^{[n-k]}\right).$$
	The inclusions are denoted by 
	$$j_0: Z_0^{[n]} \inclu \cZ \text{ and }$$
	$$j_1^{[k,n]}: Z_1^{[k]}\times (X_2/D)^{[n-k]} \inclu \cZ, \ \ 
	j_2^{[k,n]}: (X_1/D)^{[k]}\times Z_2^{[n-k]} \inclu \cZ,$$
	
	If $\cL$ is a line bundle on $\cX_U$, then define $\cL^{[n]}=Q_*P^* \cL$. 
	\item The projections from universal closed subscheme are
	$$\xymatrix{
Z_i^{[k]}  \ar[d]^{p_1^{[k]}} \ar[r]^(.4){q_i^{[k]}} & (X_i/D)^{[k]}   \\
 X_i/D&  }  
 \text{ and } 
 \xymatrix{
Z_0^{[n]}  \ar[d]^{p_0^{[n]}} \ar[r]^{q_0^{[n]}} & X_0^{[n]}  \\
 X_0 &  }   $$  
 and for line bundles $L_i$ on $X_i$, define $$L_i^{[k]}=(q_i^{[k]})_*(p_i^{[k]})^* L_i\ , i=0,1,2.$$ 

 \item $\pi_1^{[k,n]}$,  $\pi_2^{[k,n]}$ are the projections
 $$\xymatrix{
 (X_1/D)^{[k]}\times (X_2/D)^{[n-k]}   \ar[d]^{\pi_1^{[k,n]}} \ar[rr]^(.6){\pi_2^{[k, n]}}& & (X_2/D)^{[n-k]}  \\
 (X_1/D)^{[k]} &&  }   $$  
   
\end{enumerate}
  \end{depr}

\subsection{Lemmas}
Next, we will present several lemmas before proving Theorem \ref{dgfl}. 

\begin{lemm}\label{fguv} Suppose we have a fibered diagram 
$$
\xymatrix{
X \times_Z Y \ar[d]^{v} \ar[r]^{u}
& X \ar[d]^g \\
 Y \ar[r]^f & Z } $$
and a vector bundle $\cF$ on $X$. Suppose $g$ is finite and $f$ is a closed immersion, then
$$v_*u^*\cF\cong f^*g_*\cF.$$ 
\begin{proof} Since the question is local on $Z$ and $g$ is affine, we can assume $Z=\text{Spec} A$, $X=\text{Spec} B$ and $Y=\text{Spec} A/I$. Then this lemma follows from simple algebra 
$$M \otimes_{B}  \left(B \otimes_A A/I\right) \cong M \otimes_{A} A/I$$
as $A/I$-modules. 
\end{proof}

\end{lemm}
\begin{lemm}The sheaves defined in \ref{longdef} satisfy the following properties:
\begin{enumerate}
	\item $\Ln$ is a vector bundle of rank $3n$ on $\cXn$, 
	\item $L_i^{[k]}$ is a vector bundle of rank $3k$ on $X_i^{[k]}$, for $i=0,1,2$.
\end{enumerate}
\end{lemm}

\begin{proof}

All the proofs are similar. The morphisms $Q$, $q_n$ and $q_i^{[k]}$ are all finite and flat for $i=1,2$. 
Therefore the functors $R^1Q_*$, $R^1(q_0^{[n]})_*$ and $R^1(q_i^{[k]})_*$ are zero and the lemma follows from the  cohomology and base change theorem.  
\end{proof}

\begin{lemm} Recall that $i_0^{[n]}: X_0^{[n]} \inclu \cXn$ and $ i^{[k,n]}: (X_1/D)^{[k]} \times (X_2/D)^{[n-k]} \inclu \cXn$ are the inclusions of the fibers. Then 
\begin{enumerate}
	\item $(i_0^{[n]})^* \Ln= L_0^{[n]},$
	\item $(i^{[k,n]})^* \Ln= (\pi_1^{[k,n]})^* L_1^{[k]} \oplus (\pi_2^{[k,n]})^* L_2^{[n-k]} .$
\end{enumerate}
\end{lemm}
\begin{proof}

(1) We have the commutative diagram
  $$\xymatrix{
X_0^{[n]} \ar@{^{(}->}[r]^{i_0^{[n]}}   & \cXn  \\
Z_0^{[n]} \ar[u]^{q_0^{[n]}} \ar[d]_{p_0^{[n]}} \ar@{^{(}->}[r]^{j_0}& \cZ  \ar[u]_{Q}\ar[d]^{P} \\
X_0 \ar@{^{(}->}[r]^{i_0}   & \cX_U
} $$
Because diagram on the top satisfies the assumption of lemma \ref{fguv}, 
$$(i_0^{[n]})^*Q_*(P^*\cL)\cong (q_0^{[n]})_* (j_0)^*(P^*\cL).$$
Since $P\circ j_0 = i_0\circ p_0^{[n]}$, $j_0^*P^* \cL \cong (p_0^{[n]})^*i_0^*\cL $. By definition $Q_*P^* \cL=\cL^{[n]}$, therefore 
$$(i_0^{[n]})^* \cL^{[n]} \cong (q_0^{[n]})_* (p_0^{[n]})^*i_0^* \cL\cong (q_0^{[n]})_* ({p_0^{[n]}})^* L_0 \cong L_0^{[n]}. $$	

(2) The following diagram is commutative: 
$$	\xymatrix{
(X_1/D)^{[k]}  \  && \ar[ll]_<<<<<<{\pi_1^{[k,n]}} (X_1/D)^{[k]} \times (X_2/D)^{[n-k]} \ \ar@{^{(}->}[rr]^>>>>>>>>>>{i^{[k,n]}}&& \cX^{[n]}  \\ \\
Z_1^{[k]} \ar[uu]^{q_1^{[k]}} \ar[dd]_{p_1^{[k]}}&& 
\ar[ll]^{\mathrm{pr_1}} Z_1^{[k]}\times (X_2/D)^{[n-k]}  \ar[rr]^>>>>>>>>>>>>>{j_1^{[k,n]}}\ar[dd]_{\phi_1^{[k]}=\alpha_1\circ p_1^{[k]}\circ \mathrm{pr_1}} \ar[uu]^{q_1^{[k]}\times id}
\ar[dd]
&&\cZ  \ar[uu]_{Q}\ar[dd]^{P}\\ \\
X_1/D \ar@{^{(}->}[rr]^{\alpha_1}   && X_1\cup_D X_2\ar@{^{(}->}[rr]^{i_{\infty}}   && \cX_U
}$$

Since $P \circ (j_1^{[k,n]}) = i_{\infty}\circ \phi_1^{[k]}$, 
$(j_1^{[k,n]})^*(P^*\cL) =(\phi_1^{[k]})^*  (i_{\infty}^*\cL)$. Furthermore, apply Lemma \ref{fguv} on the following diagram 
$$\xymatrix{
(X_1/D)^{[k]} \times (X_2/D)^{[n-k]} \ar@{^{(}->}[rrr]^(.7){i^{[k,n]}} &&& \cX^{[n]}\\ \\
Z_1\times (X_2/D)^{[n-k]} \cup (X_1/D)^{[k]}\times Z_2^{[n-k]} \ar@{^{(}->}[rrr]^(.7){j_1^{[k,n]}\cup j_2^{[n-k,n] }}
\ar[uu]^{q_1^{[k]}\times id \cup id \times q_2^{[n-k]}} &&& \cZ \ar[uu]_Q }$$
yields
$$(i^{[k,n]})^* Q_* (P^* \cL)\cong 
(q_1^{[k]}\times id \cup id \times q_2^{[n-k]})_*(j_1^{[k,n]}\cup j_2^{[n-k,n]})^*(P^* \cL).$$
By definition $Q_*P^*\cL=\cL^{[n]}$, thus
\begin{align*}
(i^{[k,n]})^*\cL^{[n]}\cong& (q_1^{[k]}\times id \cup id \times q_2^{[n-k]})_*\left((j_1^{[k,n]})^* (P^* \cL)\cup (j_2^{[n-k,n]})^*(P^* \cL)\right)\\
\cong&  (q_1^{[k]}\times id \cup id \times q_2^{[n-k]})_* \left((\phi_1^{[k]})^*(i_{\infty})^*\cL \cup (\phi_2^{[n-k]})^*(i_{\infty})^*\cL\right)\\
= &(q_1^{[k]}\times id)_*\left((\phi_1^{[k]})^*(i_{\infty})^*\cL\right)\oplus  (id \times q_2^{[n-k]})_* \left((\phi_2^{[n-k]})^*(i_{\infty})^*\cL\right),
\end{align*}
and
\begin{align*} (q_1^{[k]}\times id)_*\left((\phi_1^{[k]})^*(i_{\infty})^*\cL\right)\cong &
(q_1^{[k]}\times id)_*\left((\tau_1^{[k,n]})^*(p_1^{[k]})^*(\alpha_1^{[k]})^*(i_{\infty})^*\cL\right)\\
\cong & (q_1^{[k]}\times id)_*\left((\tau_1^{[k,n]})^*(p_1^{[k]})^* L_1\right)\\
\cong & (\pi_1^{[k,n]})^*(q_1^{[k]})_*(p_1^{[k]})^*L_1\cong (\pi_1^{[k,n]})^*L_1^{[k]}.\end{align*}
The last isomorphism follows from applying Lemma \ref{fguv} on the upper-left part of the diagram. 
Similarly, $$(id \times q_2^{[n-k]})_* \left((\phi_2^{[n-k]})^*(i_{\infty})^*\cL\right)\cong
(\pi_2^{[n-k,n]})^*L_2^{[n-k]}.$$
Thus we can conclude 
$$(i^{[k,n]})^* \Ln= (\pi_1^{[k,n]})^* L_1^{[k]} \oplus (\pi_2^{[k,n]})^* L_2^{[n-k]}. $$
\end{proof}
%%%%%%%%%%%%%%%%%%%%%%%%%%%%%%%
Suppose $S$ is a smooth algebraic surface and $r\in \NN$. 
Recall that in Definition \ref{S2r} we have defined  $W_0^{3r}$ and 
$W^{3r}$ as subschemes of $S^{[3r]}$. When  $S=X_0$, ($X_1/D, X_2/D$ respectively),  the corresponding closed subschemes are denoted by $W_{X_0}^{3r}$, ($W_{X_1/D}^{3r}, W_{X_2/D}^{3r}$ respectively)

\begin{lemm}
There is a family of closed subschemes $\cW^{3r} \subset \cX^{[3r]}$ such that 
$$\cW^{3r} \cap X_0^{[3r]} =W_{X_0}^{3r}$$ and 
$$\cW^{3r} \cap \left({(X_1/D)^{[m]}\times (X_2/D)^{[3r-m]}}\right)= 
\begin{cases} \text{empty set} & \text{if $m$ is not divisible by 3},\\
               W_{X_1/D}^{3k}\times W_{X_2/D}^{3(r-k)} &  \text{if $m=3k$, $k\in \NN$} \end{cases}$$
Furthermore, $\cW^{3r}$ is flat over $U$ via the composition  $\cW^{3r} \inclu \cX^{[3r]}\to U$.
\end{lemm}

\begin{proof} 
Let $\cV$ be the set of points $\displaystyle \coprod_{i=1}^{r}\text{Spec}(\cO_{S, x_i}/m^2_{S, x_i})$ in $\cX^{[3r]}$ where $S$ is a smooth fiber of $\cX$  and $x_1, \ldots , x_r$ are distinct closed points on $S$. Define $\cW^{3r} $ to be the closure of $\cV$ in $\cX^{[3r]}$. 
Since limit points of family in $\cV$ are disjoint zero-dimensional subschemes of length a multiple of three, there is neither limit point nor interior point on ${(X_1/D)^{[m]}\times (X_2/D)^{[3r-m]}}$  if $m$ is not divisible by 3. 

Suppose $b$ is a closed point in $\cW^{3r} $, then $b$ is the limit of points in $\cV$. 
As a closed subscheme of $\cX_U[n]$, $b$ is the disjoint union of zero-dimensional closed subschemes $b_i$ and every $b_i$ is supported on a single point. If we choose small disjoint open neighborhoods of $b_i$, points in $\cV$ approaching $b$ in those neighborhoods are disjoint families of points in $\cV$ approaching $b_i$ for every $i$.  Therefore it suffices to assume $b$ is supported on a single point $p$ and prove $b$ is in $W_{X_0}^{3r}$ or $W_{X_1/D}^{3k}\times W_{X_2/D}^{3(r-k)}$.

Suppose $b$ is a closed point in $\cW^{3r} \cap X_0^{[3r]}$, then $p$ is a closed point of  $X_0$. The family $\cX_U\to  U$  is locally trivial near $p$. 
Thus there exists an analytic neighborhood of $p$ in $\cX_U$ which is a product of $U_1$ and $U_2$, such that $U_1$ is open in $X_0$ and isomorphic to an open set of $\CC^2$ and $U_2$ is open in $U$. 
Suppose $v_i$ is a family of closed points of $\cV$ with limit $b$, we can assume $v_i$ is sufficiently close to $b$ and belongs to $U_1\times U_2$.  
Define $v_i'$ to be the image of $v_i$ under the first projection $\mathrm{pr_1}: U_1\times U_2 \to U_1$. 
$v_i'$ is a disjoint union of $\text{Spec}(\cO_{U_1, \mathrm{pr_1}(x_i)}/m^2_{U_1, \mathrm{pr_1}(x_i)})$. 
In addition, $v_i'$ is in $U_1\subset X_0$ and the limit is still $b$. 
Therefore we can conclude that $b$ is in $W_{X_0}^{3r}$. 

In the case $b$ is on $\cW^{3r} \cap (X_1/D)^{[3k]}$ for some $k$, the zero-dimensional closed subscheme $b$ belongs to a component of $(X_1)[m]_0$ for some $m\in \NN$. Call the component $S$, 
by construction of the moduli stack the support $p$ is a nonsingular point of $S$.  
Since $S$ is contained in the fiber of $\infty$ in $\cX_U[m]$ and the family $\cX_U[m] \to U$ is also locally trivial near $p$. 
We can use the projection argument again and show $b$ is the limit of points of the form $\coprod_{i=1}^{k}\text{Spec}(\cO_{S, y_i}/m^2_{S, y_i})$ where $y_1, \ldots , y_k$ are distinct points in $S$. 
This shows $b$ belongs to $W_{X_1/D}^{3k}$. The case  $b\in \cW^{3r} \cap (X_2/D)^{[3(r-k)]}$ can be proved similarly. 

The previous argument shows that
$$ \cW^{3r} \cap X_0^{[3r]} \subset W_{X_0}^{3r} ,$$ and 
\begin{align*}
\cW^{3r} \cap \left({(X_1/D)^{[m]}\times (X_2/D)^{[3r-m]}}\right)&= \phi & \text{if $m$ is not divisible by 3},\\
\cW^{3r} \cap \left({(X_1/D)^{[m]}\times (X_2/D)^{[3r-m]}}\right)&\subset 
               W_{X_1/D}^{3k}\times W_{X_2/D}^{3(r-k)} &  \text{if $m=3k$, $k\in \NN$}. \end{align*}

On the other hand, $W_{X_0}^{3r}$ is the closure of $\cV\cap X_0^{[3r]}$. This is because $\cW^{3r}$ contains $\cV\cap X_0^{[3r]}$ and is a closed set, therefore $\cW^{3r}$ must contain $W_{X_0}^{3r}$. 
Similarly, $W_{X_1/D}^{3k}\times W_{X_2/D}^{3(r-k)}$ is the closure of $\cV \cap \left((X_1/D)^{[3k]}\times (X_2/D)^{[3(r-k)]}\right)$. Because   $\cW^{3r}$ contains $\cV \cap \left((X_1/D)^{[3k]}\times (X_2/D)^{[3(r-k)]}\right)$ and is closed, $W_{X_1/D}^{3k}\times W_{X_2/D}^{3(r-k)}$ is a subset of $\cW^{3r}$. 
This finished the first part of the proof. 

Since every point in $\cV$ is a reduced closed point in $\cX^{[3r]}$ and $\cV$ is dense in $\cW^{[3r]}$, $\cW^{[3r]}$ is reduced. Second, for every point $b$ in $\cW^{[3r]}$  there is a section $U\to \cW^{[3r]}$ which passes through $b$. Therefore $\cW^{[3r]}$ is irreducible and dominates $U$ which implies $\cW^{[3r]}\to U$ is flat.

\end{proof}  

\subsection{Relative generating functions}
\begin{defn}
Define 
$$d_0(X_i/D, L_i) =1,\  d_k(X_i/D, L_i) =\int_{W_{Xi/D}^{3k}} c_{2k}(L_i^{[3k]}) \text{ for } k\geq 1  $$   
$$\text{ and } \phi(X_i/D, L_i)(x)=\sum_{k=0}^{\infty}  d_k(X_i/D, L_i) x^k $$
to be the relative enumerative number and the relative generating function of $[X_i/D,L_i]$ for $i=1, 2$. 
\end{defn}

\begin{prop} \label{012} Suppose $[X_0, L_1]-[X_1, L_1]-[X_2, L_2]+[X_3, L_3]$ is an extended double point relation, then
$$\phi(X_0, L_0)= \phi(X_1/D, L_1)\cdot \phi(X_2/D, L_2).$$
\end{prop}

\begin{proof}
For $r\in \NN$, recall $X_0^{[3r]}$ is the fiber of $\cX^{[3r]}$ over $0$ and 
$$\cL^{[3r]}|_{X_0^{[3r]}}=L_0^{[3r]}, \  \cW^{3r} \cap X_0^{[3r]}=W_{X_0}^{3r}.$$ 
The fiber of $\cX^{[3r]}$ over $\infty$ is $$\cX^{[3r]}|_{\infty}=\cup_{m=0}^{3r} (X_1/D)^{[m]} \times (X_2/D)^{[3r-m]}\,,$$
$$\cL^{[3r]}|_{(X_1/D)^{[m]} \times (X_2/D)^{[3r-m]}} =
(\pi_1^{[m,3r]})^* L_1^{[m]} \oplus (\pi_2^{[m,3r]})^* L_2^{[3r-m]},$$
$$\cW^{3r} \cap \cX^{[3r]}|_{\infty}=\cup_{k=0}^{3r} W_{X_1/D}^{3k} \times W_{X_2/D}^{3(r-k)}.$$ 
Recall that $U$ is open in $\PP^1$, $\cX^{[3r]} \to U$ is a flat family and $\cW^{3r} \to U$ is flat. Thus $[\cW^{3r}|_{0}]$ and $[\cW^{3r}_{\infty}]$, the classes of fibers of $\cW^{3r}$ over $0$ and $\infty$, satisfy
\begin{align*}c_{2r}(\cL^{[3r]})[\cW^{3r} |_{0}]=c_{2r}(\cL^{[3r]})[\cW^{3r} |_{\infty}].\end{align*}
\begin{align*}c_{2r}(\cL^{[3r]})[\cW^{3r} |_{0}]=&c_{2r}(\cL^{[3r]})[W_{X_0}^{3r}]
=\int_{W_{X_0}^{3r}} c_{2r}(L_0^{[3r]})=d_r(X_0, L_0),
\end{align*}
\begin{align*}&c_{2r}(\cL^{[3r]})[\cW^{3r} |_{\infty}]
=c_{2r}(\cL^{[3r]})[\cup_{k=0}^{3r} W_{X_1/D}^{3k} \times W_{X_2/D}^{3(r-k)}]\\
=&\sum_{k=0}^{r} c_{2r}((i^{[3k, 3r]})^*\cL^{[3r]})[ W_{X_1/D}^{3k} \times W_{X_2/D}^{3(r-k)}]\\
=&\sum_{k=0}^{r} c_{2r}\left((\pi_1^{[3k,3r]})^* L_1^{[3k]} \oplus (\pi_2^{[3k,3r]})^* L_2^{[3(r-k)]}\right)
[W_{X_1/D}^{3k} \times W_{X_2/D}^{3(r-k)}]\\
= &\sum_{k=0}^{r}\left(\sum_{l=0}^{2r} c_l\left((\pi_1^{[3k,3r]})^* L_1^{[3k]}\right)c_{2r-l}\left((\pi_2^{[3k,3r]})^* L_2^{[3(r-k)]}\right)\right)[W_{X_1/D}^{3k} \times W_{X_2/D}^{3(r-k)}].
\end{align*}
For dimension reasons, the only nonzero terms are:
\begin{align*}
&\sum_{k=0}^{r} \left(c_{2k}\left((\pi_1^{[3k,3r]})^* L_1^{[3k]}\right) c_{2(r-k)}\left((\pi_2^{[3k,3r]})^* L_2^{[3(r-k)]}\right)\right)
[W_{X_1/D}^{3k} \times W_{X_2/D}^{3(r-k)}]\\
=&\sum_{k=0}^{r}  \int_{W_{X_1/D}^{3k}} c_{2k}(L_1^{[3k]}) \int_{W_{X_2/D}^{3(r-k)}} c_{2(r-k)}(L_2^{[3(r-k)]})\\
=&\sum_{k=0}^{r}  d_k(X_1/D, L_1)\cdot d_{r-k}(X_2/D, L_2).
\end{align*}
Therefore for every positive integer $r$, 
$$d_r(X_0, L_0)=\sum_{k=0}^{r}  d_k(X_1/D, L_1)\cdot d_{r-k}(X_2/D, L_2),$$
$$\phi(X_0, L_0)=\phi(X_1/D, L_1)\cdot \phi(X_2/D, L_2).$$ 

\end{proof}\smallskip

To obtain a formula for absolute generating function, the relative surfaces can be closed up by adding projective bundles.

\begin{coro}\label{closeup} Let $[X, L]$ be a pair in $\SL$. Suppose $D$ is a smooth curve in $X$, $N$ is the normal bundle of $D$ in $X$ and $L\cdot D$ is the intersection number of $c_1(L)$ and $D$, then

$$\phi(X,L)= \phi(X/D, L) \cdot \phi(\PP_N/D, (L\cdot D)f).$$

\end{coro}\smallskip

\begin{proof}Let $\cX$ be the blow up of $X\times \PP^1$ along $D \times  \{\infty\}$, then $\cX$ is the deformation to the normal cone. Let $\pi: \cX \to X\times \PP^1 \to \PP^1$ be the composition of blow-up morphism and projection to $\PP^1$. $\cX \to \PP^1$ is a smooth family of surfaces with general fiber $X$ and $\pi^{-1}(\infty)= X \cup_D \PP_N$.  Apply Proposition \ref{012} to $\cX$, then we can obtain the desired formula.

\end{proof}

For relative surfaces $X_1/D$ and $X_2/D$, if we let the normal bundle of $D$ in $X_1$ be $N$, then the normal bundle of $D$ in $X_2$ is the dual bundle $N^{\vee}$. 
Recall $$\PP_{N}= \PP(\cO_D \oplus N ), \ \PP_{N^{\vee}}=\PP(\cO_D \oplus  N^{\vee}). $$ 
We fix the embedded curves 
$$D_0=\PP(N) \subset \PP_{N}, \ D_{\infty}= \PP(\cO_D) \subset \PP_{N}, \ 
D_0^{\vee}=\PP(N^{\vee}) \subset \PP_{N^{\vee}}, \ D_{\infty}^{\vee}= \PP(\cO_D) \subset \PP_{N^{\vee}}.$$

The curve $D$ can be naturally embedded to the four curves in $\PP_{N}$ and $\PP_{N^{\vee}}$. 
In addition, because $L_1$ and $L_2$ are the restriction of line bundle $\cL$ on $\cX$,

$$L\cdot D_1=c_1(L_1) \cap [X_1]\cap [X_2]=c_1(\cL)\cap [D] \cap [X_1]\cap[ X_2]=L\cdot D_2.$$ 
This implies that in Definition \ref{dp}, if $f$ is the fiber class in $X_3\to D$, the line bundle $L_3$ equals $(L_1\cdot D)f=(L_2\cdot D)f=(c_1(\cL)\cap D)f$. Therefore we will also use $L_3$ to denote the line bundle  $(L_1\cdot D)f$ on $\PP_N$ and $(L_2\cdot D)f$ on $\PP_{N^\vee}$.

\begin{coro} Let $X_i$, $D$, $\PP_{N}$, $\PP_{N^{\vee}}$ be as above and in Definition \ref{dp}. The generating functions satisfy
$$\phi(X_1, L_1)= \phi(X_1/D, L_1)\cdot \phi(\PP_{N}/D_0, L_3), $$
$$\phi(X_2, L_2)= \phi(X_2/D, L_2) \cdot \phi(\PP_{N^{\vee}}/D_0^{\vee}, L_3)$$
\end{coro}

\begin{proof} Apply Corollary \ref{closeup} to $X_1/D$ and $X_2/D$. \end{proof}

Let $X_0= \PP_{N}$, $D=D_0$ and $\cX$ be the the blowup of $\PP_{N}\times \PP^1$ at $D \times \{0\}$. 
Then the general fiber is $\PP_{N}$ and $\pi^{-1}(\infty)= \PP_{N}\cup_D \PP_{N}$ where $D$ embeds into the first $\PP_{N}$ as $D_0$ and embeds into the second $\PP_{N}$ as $D_{\infty}$. Let  $\cL$ be the pullback of $L_3$ from the composition $\cX\to \PP_{N}\times \PP^1 \to \PP_{N}$.  Proposition \ref{012} implies 

\begin{coro}
$$\phi(\PP_{N}, L_3)= \phi(\PP_{N}/D_0, L_3) \cdot 
\phi(\PP_{N}/D_{\infty}, L_3).$$
\end{coro}\smallskip

Now we are ready to prove Theorem \ref{dgfl}

\begin{proof}[Proof of Theorem \ref{dgfl}]
From Proposition \ref{012} and its corollaries:

\begin{align*}\phi(X_0, L_0)&= \phi(X_1/D, L_1) \cdot \phi(X_2/D, L_2),\\
\phi(X_1/D, L_1) \cdot \phi(\PP_{N}/D_0, L_3) &=\phi(X_1, L_1),\\
\phi(X_2/D, L_2) \cdot  \phi(\PP_{N^{\vee}}/D_0^{\vee}, L_3)&=\phi(X_2, L_2),\\
\phi(\PP_{N}, L_3)&= \phi(\PP_{N}/D_0, L_3) \cdot 
\phi(\PP_{N}/D_{\infty}, L_3).\end{align*}

There is a canonical isomorphism between $\PP_N$ and $\PP_{N^{\vee}}$ which maps $D_{\infty}\subset\PP_N $ to $D_0\subset \PP_{N^{\vee}}$ and keeps the line bundle $L_3$ unchanged. 
This implies $$\phi(\PP_{N^{\vee}}/D_0^{\vee}, L_3)=\phi(\PP_{N}/D_{\infty}, L_3).$$ 
Then the theorem is proved by multiplying all equations. 

\end{proof}
  
\section{Universality Theorems and Generating functions}\label{chap4}

\subsection{Outline}
In this Section, we will prove Theorem \ref{univpoly}, Theorem \ref{DG2} and Theorem \ref{power} by combining the degeneration formula (Theorem \ref{dgfl}) with the structure of algebraic cobordism group $\SL$ (Theorem \ref{slfl}). 

Recall that for any smooth projective surface $S$ and line bundle $L$ on $S$, we defined and studied the enumerative number $d_r(S,L)$ and generating function 
$$\phi(S,L)(x)=\sum_{r=0}^{\infty} d_r(S,L)\, x^r$$
in Section \ref{degformula}. By Theorem \ref{dgfl}, this function $\phi$ induces a homomorphism from $\SL$ to $(\qqx, \ \cdot \ )$. 
Theorem \ref{slfl} proves that $\SL$ is four-dimensional and the only invariants are \fourtopand. 
Combining these two results, we show that $\phi(S,L)(x)$ only depends on these four topological numbers and has a multiplicative structure. 
 
\subsection{Proof of Theorem 1.3 and Theorem 1.1}

\begin{prop}\label{thm:4.1} There exist four series $A_1$, $A_2$, $A_3$ and $A_4$ in $\qqx$ such that 
$$\phi(S,L)(x)=A_1^{L^2}A_2^{LK_S}A_3^{c_1(S)^2}A_4^{c_2(S)}. $$ More explicitly,   
\begin{align*}
A_1(x)&= \phi(\PP^2, \cO)^{-1} \phi(\PP^2, \cO(1)) \phi(\ptp, \cO)^{\frac32 } \phi(\ptp, \cO(1,0))^{-\frac32 },\\    
A_2(x)&= \phi(\ptp, \cO)^{\frac12 } \phi(\ptp, \cO(1,0))^{-\frac12},\\
A_3(x)&= \phi(\PP^2, \cO)^{-\frac13 } \phi(\ptp, \cO)^{-\frac14},\\
A_4(x)&= \phi(\PP^2, \cO)^{-\frac23}\phi(\ptp, \cO)^{\frac34}.
\end{align*} 
\end{prop}

\begin{proof}By Theorem \ref{slfl}, the class of $[S,L]$ in $\SL$ is
\begin{align}\label{eq:4.1}[S,L]= a_1[\PP^2, \cO]+a_2[\PP^2, \cO(1)]+a_3 [\ptp, \cO]+a_4[\ptp, \cO(1,0)] \end{align}
where 
\begin{align*}
a_1&= -L^2+\frac{c_1(S)^2+c_2(S)}{3}-c_2(S), &a_2&= L^2,\\
a_3&=L^2+\frac{LK+L^2}{2}-\frac{c_1(S)^2+c_2(S)}{4}+c_2(S), & a_4&=-L^2-\frac{LK+L^2}{2}.\end{align*}

Since the generating function $\phi$ is a homomorphism from $\SL$ to $(\qqx, \ \cdot \ )$ (Theorem \ref{dgfl}). We obtain
\begin{align*} 
\phi(S,L)(x) =\phi(\PP^2, \cO)^{a_1}\phi(\PP^2, \cO(1))^{a_2}\phi(\ptp, \cO)^{a_3}\phi(\ptp, \cO(1,0))^{a_4}. \end{align*}
Plugging in the values of $a_i$ and grouping the functions according to \fourtopand\ complete the proof. 
\end{proof}

Another way to prove Proposition \ref{thm:4.1} is using the following diagram 

$$   \xymatrix{
& \qqx& \\
&&\\
\SL \  \ar@{=}[rr]_{(L^2, \,LK, \, c_1(S)^2, \,c_2(S))} \ar[ruu]^{\phi}&&
\QQ^4 \ar@{-->}[luu]}
$$

Since $\phi$ is a homomorphism and $\SL \to \QQ^4$ is an isomorphism, they induce a homomorphism from $\QQ^4$ to $\qqx$. This implies $\phi(S,L)(x)$ is an universal power series in \fourtopand. Furthermore, let the the image of the standard basis $\mathbf{e_i}$ of $\QQ^4$ be $A_i$, then the commutative diagram implies 
$$\phi(S,L)(x)=A_1^{L^2}A_2^{LK_S}A_3^{c_1(S)^2}A_4^{c_2(S)}. $$

\begin{rem}One can see the same proof is still valid if we replace 
$$\{[\PP^2, \cO], [\PP^2, \cO(1)], [\ptp, \cO], [\ptp, \cO(1,0)]\}.$$
by any basis of $\SL$. 
Thus the four series $A_i's$ can also be expressed in the generating function $\phi$ of any basis.\end{rem}

\begin{coro}The enumerative number $d_r(S,L)$ is a universal polynomial of degree $r$ in \fourtopand\ for all projective smooth surfaces $S$ and line bundles $L$ on $S$.  \end{coro}
\begin{proof}
For all pairs $[S,L]$, $d_r(S,L)$ is the coefficient of $x^r$ in 
$$\phi(S,L)(x)= A_1^{L^2}A_2^{LK_S}A_3^{c^2_1(S)}A_4^{c_2(S)}\,.$$ 
After expanding the series, one can see that the coefficient of $x^r$ in $\phi(S,L)(x)$ is a universal polynomial of $(L^2, LK,  c_1(S)^2, c_2(S))$ of degree $r$. 
\end{proof}

Next, we prove \Gott's conjecture: 

\begin{univpoly}
For every integer $r\geq 0$, there exists a universal polynomial $T_r(x,y,z,t)$ of degree $r$ with the following property: given a smooth projective surface $S$ and a $(5r-1)$-very ample (5-very ample if $r= 1$) line bundle $L$ on $S$, a general $r$-dimensional sublinear system of $|L|$ contains exactly $T_r(L^2, LK, c_1(S)^2, c_2(S))$ $r$-nodal curves.  
\end{univpoly}

\begin{proof}
$d_r(S,L)$ is always a universal polynomial of degree $r$ for all line bundles $L$ on $S$. 
If $L$ is $(5r-1)$-very ample, Proposition \ref{Gott'sthm} shows $d_r(S, L)$ equals the number of $r$-nodal curves in $[S,L]$. Thus the number of $r$-nodal curves in $[S,L]$ is given by a universal polynomial of degree $r$ when $L$ is $(5r-1)$-very ample. 
\end{proof}

\begin{coro}\label{coro:4.3}The generating functions $\phi(S,L)(x)$ and $T(S,L)(x)$ (defined in Section \ref{GYZ} as $\sum_{r=0}^{\infty} T_r(L^2, LK, c_1(S)^2, c_2(S))x^r$) are equal. 
\end{coro}

\begin{power}There exist universal power series $A_1$, $A_2$, $A_3$, $A_4$  in $\qqx$ such that the generating function $T(S,L)$ has the form
$$T(S,L)=A_1^{L^2}A_2^{LK_S}A_3^{c_1(S)^2}A_4^{c_2(S)}.$$\end{power}

\begin{proof}This follows from Proposition \ref{thm:4.1} and Corollary \ref{coro:4.3}. \end{proof}

\subsection{Proof of the \Gott-Yau-Zaslow Formula}
Let $\tau$ be on the complex upper half plane, $q=e^{2\pi i\tau}$. 
The Eisenstein series $$G_2(\tau):=\frac{-1}{24}+\sum_{n>0}\left(\sum_{d|n} d\right)q^n$$ 
is a quasi-modular form. Denote by $D$ the differential operator $D :=\frac{1}{2\pi i}\frac{d}{d\tau}=q\frac{d}{dq}$. 
Since the ring of quasimodular forms is closed under differentiation, $DG_2$ and $D^2G_2$ are also quasi-modular. 
In addition, let $$\Delta(\tau)=q\prod_{k>0}(1-q^k)^{24}=\eta(\tau)^{24}$$ where 
$\eta(\tau)$ is the Dedekind $\eta$ function,  $\Delta(\tau)$ is a modular form.  
We define the generating function 
$$\gamma(S,L)(q)=\sum_{r=0}^{\infty}T_r(L^2, LK, c_1(S)^2, c_2(S))(DG_2(\tau))^r. $$ 
An equivalent form of Bryan and Leung's result (\cite{BL}, Theorem 1.1) states that for a generic K3 surface $S$ and a primitive class $L$, 
$$\gamma(S,L)(q)= \frac{(DG_2(\tau)/q)^{\chi(L)}}{\Delta(\tau)D^2G_2(\tau)/q^2}. $$ 
More generally, the generating function $\gamma(S,L)(q)$ for all algebraic surfaces is given by
  
\begin{GYZ}[The \Gott-Yau-Zaslow Formula]
There exist universal power series $B_1$, $B_2$ in $q$ such that
$$\sum_{r=0}^{\infty}T_r(L^2, LK, c_1(S)^2, c_2(S))(DG_2(\tau))^r=
\frac{(DG_2(\tau)/q)^{\chi(L)}B_1(q)^{K_S^2}B_2(q)^{LK_S}}{(\Delta(\tau)D^2G_2(\tau)/q^2)^{\chi(\cO_S)/2}}.$$
\end{GYZ} 

\begin{proof}
Since $(K_S^2, LK_S,  \chi(L), \chi(\cO_S))$ and $(\fourtop)$ determine each other by linear relations, $\phi(S,L)$ is multiplicative in $(K_S^2, LK_S, \chi(L), \chi(\cO_S))$ and  so as $\gamma(S,L)(q)= \phi(S,L)(DG_2)$. Therefore we can let $$\gamma(S,L)(q)=B_1(q)^{K_S^2}B_2(q)^{LK_S}B_3(q)^{\chi(L)}B_4(q)^{\chi(\cO_S)}.$$

For all generic K3 surfaces $S$ and primitive classes $L$,
$$\gamma(S,L)(q)= \frac{(DG_2(\tau)/q)^{\chi(L)}}{\Delta(\tau)D^2G_2(\tau)/q^2} .$$ 
These pairs can achieve infinitely many distinct $\chi(L)$, thus $$B_3(q)=DG_2(\tau)/q \ \ \text{ and } B_4(q)=(\Delta(\tau)D^2G_2(\tau)/q^2)^{1/2}$$ and the result follows.  
\end{proof}

\begin{rem}Although in \cite{BL} they found the formula for generating function of  primitive classes on K3 surfaces with Picard number one, it is not trivial to see why the coefficients equals our universal polynomials. 
For each $r$, we have to show that there exist infinitely many K3 surfaces with Picard number one and primitive classes $L$ such that $L$ is $(5r-1)$-very ample and also have distinct intersection numbers $L^2$. It can be proved by using Lemma 5.3 in \cite{KP01}, which shows that on a K3 surface with Picard number one, a primitive class $L$ is $(5r-1)$-very ample if $L^2> 20r$.  
Therefore Bryan and Leung's generating function is exactly $\phi(S,L)(x)$. 
\end{rem}

There is a reformulation of Theorem \ref{GYZ} given by G\"{o}ttsche (\cite{Gott}, Remark 2.6), which is the original version of Bryan and Leung's formula. 

\begin{coro} For all $l$, $m$, $r$, define
$$n_r^S(l, m)=T_{l+\chi(\cO_S)-1-r}^S (2l+m, m). $$ Then 

\begin{align*}
\sum_{l\in \ZZ} 	n_r^S(l, m) q^l= B_1(q)^{K_S^2}B_2(q)^m (DG_2(\tau))^r\frac{D^2G_2(\tau)}{(\Delta(\tau)D^2 G_2(\tau)/q^2)^{\chi(\cO_S)/2}}.
\end{align*}
\end{coro}
If $L$ is sufficiently ample with respect to $\delta=\chi(L)-1-r$, then $n^S_r((L^2-LK_S)/2, LK_S)$ counts the $\delta$-nodal curves in a general $r$-codimensional sublinear system of $|L|$. This reformulation provides another (possibly better) way to find the closed form of $B_1(q)$ and $B_2(q)$.

\end{document}